\documentclass[12pt]{amsart}

\usepackage{frontmatter_mim_nanopteron}

\begin{document}

\title[Small Mass Nanopterons in MiM Lattices]{Small Mass Nanopteron Traveling Waves in Mass-in-Mass Lattices with Cubic FPUT Potential}

\author{Timothy E. Faver}
\address{Mathematical Institute, University of Leiden, P.O. Box 9512, 2300 RA Leiden, The Netherlands, {\tt{t.e.faver@math.leidenuniv.nl}}}

\keywords{FPU, FPUT, granular chain, granular metamaterial, mass-in-mass lattice, mass-with-mass lattice, nanopteron, solitary wave}

\thanks{I thank Doug Wright for originally suggesting this problem to me and for his helpful comments and perspectives, as well as Atanas Stefanov for his interest and encouragement.
I acknowledge support from the Netherlands Organization for Scientific Research (NWO) (grant 639.032.612).}

\subjclass[2010]{Primary 35C07, 37K60; Secondary 35B20, 35B40}

\maketitle

\begin{abstract}
The mass-in-mass (MiM) lattice consists of an infinite chain of identical beads that are both nonlinearly coupled to their nearest neighbors and linearly coupled to a distinct resonator particle; it serves as a prototypical model of wave propagation in granular crystals and metamaterials.
We study traveling waves in an MiM lattice whose bead interaction is governed by the cubic Fermi-Pasta-Ulam-Tsingou potential and whose resonator mass is small compared to the bead mass.
Excluding a countable number of ``antiresonance'' resonator masses accumulating at 0, we prove the existence of nanopteron traveling waves in this small mass limit.
The profiles of these waves consist of the superposition of an exponentially localized core and a small amplitude periodic oscillation that itself is a traveling wave profile for the lattice.
Our arguments use functional analytic techniques originally developed by Beale for a capillary-gravity water wave problem and recently employed in a number of related nanopteron constructions in diatomic Fermi-Pasta-Ulam-Tsingou lattices.
\end{abstract}

\section{Introduction}

\subsection{The mass-in-mass lattice}
Consider an infinite chain of particles, which we call ``beads,'' each normalized to have mass 1, arranged on a horizontal line and connected by identical springs whose potential is $\V \in \Cal^1(\R)$.
Linearly coupled to each bead is a ``resonator'' particle of mass $\mu > 0$.
The springs connecting the beads to the resonators are all identical and exert the force $\kappa{r}$ when stretched a distance $r$.
This construct is called a mass-in-mass (MiM)\footnote{Some applications specify that the resonator is coupled to the bead internally, in which case the resonator is genuinely ``in'' the bead; others maintain the resonator outside the bead and call their lattice a ``mass-with-mass'' lattice.} lattice.

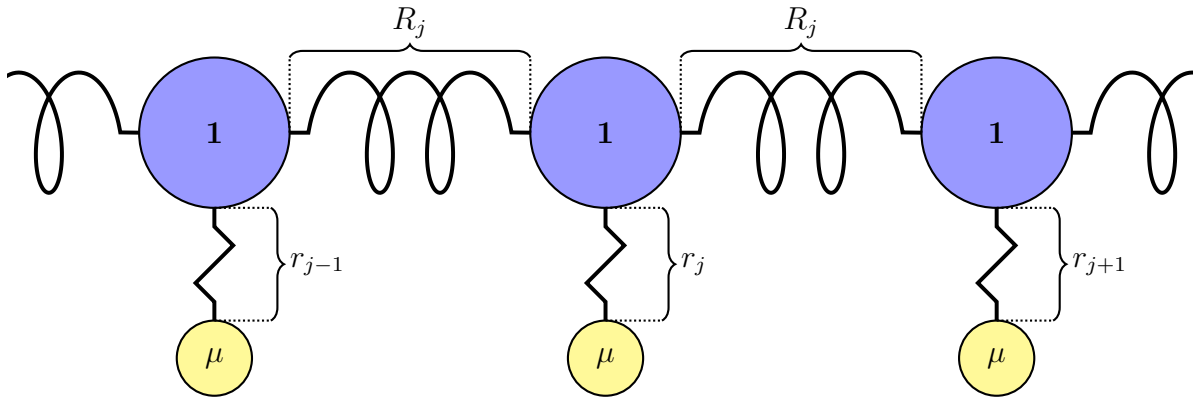
\begin{figure}
\[
\begin{tikzpicture}[thick]

\begin{scope}
\clip (-2.75,-3.55) rectangle (13,2);

\def\coilspan{2.7};
\def\edgespan{.25};
\def\beadrad{1};
\def\h{.8*\beadrad};

\def\M{.4};
\def\N{.55};
\def\E{.35};

\def\coil#1{
{\N+\M*\t+\E*sin(\t*pi r)+#1},
{\h*cos(\t*pi r)}
}

\draw[line width=1.5pt,domain={-.5:4.5},smooth,variable=\t,samples=100]
plot(\coil{-\coilspan-\edgespan-\beadrad})
--(-\beadrad,0);

\bead{0}{\beadrad};
\linearcoil{0};
\resonator{0};
\lowerlabels{0}{$r_{j-1}$};
\upperlabels{0}{$R_j$};

\draw[line width=1.5pt,domain={-.5:4.5},smooth,variable=\t,samples=100]
(\beadrad,0)--plot (\coil{\beadrad+\edgespan})
--(\beadrad+2*\edgespan +\coilspan,0);

\bead{2*\beadrad+2*\edgespan+\coilspan}{\beadrad};
\linearcoil{2*\beadrad+2*\edgespan+\coilspan};
\resonator{2*\beadrad+2*\edgespan+\coilspan};
\lowerlabels{2*\beadrad+2*\edgespan+\coilspan}{$r_j$};
\upperlabels{2*\beadrad+2*\edgespan+\coilspan}{$R_j$};

\draw[line width=1.5pt,domain={-.5:4.5},smooth,variable=\t,samples=100]
(3*\beadrad+2*\edgespan+\coilspan,0)
--plot(\coil{3*\beadrad+3*\edgespan+\coilspan})
--(3*\beadrad+4*\edgespan+2*\coilspan,0);

\bead{4*\beadrad+4*\edgespan+2*\coilspan}{\beadrad};
\linearcoil{4*\beadrad+4*\edgespan+2*\coilspan};
\resonator{4*\beadrad+4*\edgespan+2*\coilspan};
\lowerlabels{4*\beadrad+4*\edgespan+2*\coilspan}{$r_{j+1}$};

\draw[line width=1.5pt,domain={-.5:4.5},smooth,variable=\t,samples=100]
(5*\beadrad+4*\edgespan+2*\coilspan,0)
--plot(\coil{5*\beadrad+5*\edgespan+2*\coilspan})
--(5*\beadrad+6*\edgespan+3*\coilspan,0);

\end{scope}
\end{tikzpicture}
\]
\caption{The MiM lattice}
\label{fig-mim}
\end{figure}

Denote by $U_j$ the position of the center of mass of the $j$th bead and by $u_j$ the position of the center of mass of the $j$th resonator, with $j \in \Z$.
Then Newton's law implies that 
\begin{equation}\label{original lattice equations}
\begin{cases}
\ddot{U}_j = \V'(U_{j+1}-U_j)-\V'(U_j-U_{j-1})+\kappa(u_j-U_j) \\
\\
\mu\ddot{u}_j = \kappa(U_j-u_j).
\end{cases}
\end{equation}
We will rewrite \eqref{original lattice equations} in terms of relative displacements.
Set
\[
R_j := U_{j+1}-U_j 
\quadword{and}
r_j := U_j-u_j.
\]
That is, $R_j$ is the distance between the $j$th and the $(j+1)$st bead and $r_j$ is the distance between the $j$th bead and its attached resonator.
Then
\begin{equation}\label{rel disp equations}
\begin{cases}
\ddot{R}_j = \V'(R_{j+1})-2\V'(R_j)+\V'(R_{j-1})-\kappa(r_{j+1}-r_j) \\
\\
\mu\ddot{r}_j = \mu\big[\V'(R_j)-\V'(R_{j-1})\big] -\kappa(1+\mu)r_j.
\end{cases}
\end{equation}
See Figure \ref{fig-mim} for a sketch of this set-up.
We are interested in traveling waves, and so we make the ansatz
\begin{equation}\label{tw ansatz}
r_j(t)
= \begin{cases}
p_1(j-ct), \ j \text{ is odd} \\
p_2(j-ct), \ j \text{ is even},
\end{cases}
\end{equation}
where the profiles $p_1$ and $p_2$ are functions of a real variable and $c \in \R$ is the wave speed.

We will ultimately choose $\V$ to be the ``cubic'' potential
\begin{equation}\label{cubic potential}
\V(r)
= \V_{\cubic}(r)
:= \frac{r^2}{2} + \frac{r^3}{3}.
\end{equation}
We discuss the reasons for this election, as well as a survey of other relevant potentials, in Section \ref{potentials and monatomic}.

An MiM lattice is an example of a granular crystal, or granular chain, which is a lattice of closely packed particles whose behavior chiefly depends on the interactions between neighboring particles and the contact forces that they exert on each other.  
In a prototypical experiment, a finite, one-dimensional chain of particles is constructed; the chain is then struck at one end by another particle, and the resulting motion of the chain is measured \cite{cpkd}.  
Common to these experiments is an interest in the propagation of traveling waves through the granular crystal, which are waves whose form is given by a steadily translating profile  \cite{nesterenko-propagation, nesterenko}.

The appeal of the MiM lattice in modeling arises its material tunability; one can vary the spring forces from forces with a dominant linear term (the ``precompression'' case) to highly nonlinear ones (no precompression), as well as the mass values, incorporating periodic structures and models with material defects.  
For example, the articles \cite{bonanomi-et-al} and \cite{gantzounis} describe numerical and physical experiments on a finite MiM lattice consisting of a single hollow aluminum bead with an internal solid steel resonator connected to the shell by a plastic structure; polycarbonate rods connect the bead to guide plates and permit it to move. 
The resonator acts to absorb vibrations and reinforces the structural integrity of composite materials.
The MiM lattice need not be precisely a chain of particles and springs; one incarnation is the ``woodpile'' lattice of \cite{woodpile}, where the lattice ``beads'' are rods made of fused quartz, and the lattice itself assembled by stacking the rods in an orthogonal configuration.
Here, the internal vibrations in the rods serve as the internal resonator.  
As representative of the diversity of applications of MiM lattices, we mention that they arise in constructing sensors for bone elasticity \cite{bones} and ultrasonic scans \cite{Spadoni}, determining of the setting time for cement \cite{cement}, and modeling of switches and logic gates \cite{Li}.



\subsection{Reduction to monatomic lattices}\label{potentials and monatomic}
It is clear that if we take $\mu = 0$, then the small resonator in the MiM lattice vanishes, and we are left with a monatomic Fermi-Pasta-Ulam-Tsingou \cite{fput-original, dauxois} lattice with spring potential $\V$.
Depending on the spring potential selected and the wave speed regime in which one works, the monatomic FPUT lattice is known to possess a variety of solitary traveling waves, i.e., traveling waves whose profiles vanish at spatial infinity.
See Pankov \cite{pankov} for a broad overview of traveling waves in monatomic and polyatomic lattices, including, but not limited to, solitary waves, and Kevrekidis \cite{kevrekidis} for an overview of their applications (both overlapping with and distinct from those for MiM lattices discussed previously).

Our intention is to perturb from a monatomic solitary wave, in an appropriate sense, and construct traveling wave solutions for the MiM lattice when $\mu$ is small.
We make this precise in Section \ref{traveling wave problem section}.
Solitary traveling waves for the monatomic FPUT lattice have profiles $\rho$ such that putting
\[
R_j(t) = \rho(j-ct)
\]
solves
\[
\ddot{R}_j
=  \V'(R_{j+1})-2\V'(R_j)+\V'(R_{j-1}).
\]
Equivalently, $\rho$ must satisfy
\begin{equation}\label{monatomic problem}
c^2\rho'' -A\V'(\rho) = 0,
\qquad
(A\rho)(x) := \rho(x+1)-2\rho(x) + \rho(x-1).
\end{equation}
We highlight four families of bead potentials $\V$ and wave speeds $c$ for which the monatomic FPUT lattice with spring potential $\V$ is known to have solitary traveling wave solutions.

\begin{enumerate}[label={\bf(\roman*)},ref={(\roman*)}]

\item\label{fp case}
Friesecke and Pego \cite{friesecke-pego1} show that \eqref{monatomic problem} has an exponentially localized solution for $\V \in \Cal^4(\R)$ with $\V(0) = \V'(0) = 0$, $\V''(0) > 0$, and $\V'''(0) \ne 0$, with $c$ chosen in the ``near-sonic'' regime $|c| \gtrsim \sqrt{V''(0)}$.
The cubic potential $\V_{\cubic}$ from \eqref{cubic potential}, which was one of the original potentials studied by Fermi, Pasta, Ulam, and Tsingou in \cite{fput-original}, certainly satisfies these conditions, as do the Toda \cite{toda} and Lennard-Jones potentials.
The Friesecke-Pego solution turns out to be a small perturbation in $H^1$ of a KdV $\sech^2$-type solitary wave.

\item
Friesecke and Wattis \cite{friesecke-wattis} take $\V \in \Cal^2(\R)$ to be ``superquadratic'' with $r \mapsto \V(r)/r^2$ strictly increasing on either $(-\infty,0)$ or $(0,\infty)$.
The wave speed $c$ is supersonic in that $|c| > \sqrt{\V''(0)}$ but otherwise not specified.
Their methods apply to a plethora of potentials, including the cubic FPUT.

\item
Herrmann and Matthies \cite{herrmann-matthies-asymptotic, herrmann-matthies-stability, herrmann-matthies-uniqueness} construct localized solutions in the ``high-energy'' limit, i.e., with $|c|$ large, for a class of singular convex potentials that resemble the Lennard-Jones potential.

\item
The bead potential most often associated with MiM lattices in the literature is the Hertzian potential, which we now describe in some detail.
For $r \in \R$, set
\[
[r]_+
:= \begin{cases}
r, \ r \ge 0 \\
0, \ r < 0,
\end{cases}
\]
and, for $\ep_0 \ge 0$ and $p > 1$, take $\V_{\Hertz}$ to be the potential satisfying
\begin{equation}\label{Hertzian potential}
\V_{\Hertz}'(r)
:= [\ep_0+r]_+^p.
\end{equation}
When $\ep_0 > 0$, the potential models lattices experiencing a ``precompression'' force in which there is an additional force pressing on the ends of the lattice \cite{nesterenko}.
MacKay \cite{mackay}, Ji and Hong \cite{ji-hong}, English and Pego \cite{english-pego}, and Kevrekidis and Stefanov \cite{sk12, sk13} all have established the existence of solitary wave in Hertzian monatomic lattices for different values of $p$ and $\ep_0$.

Moreover, Kevrekidis, Stefanov, and Xu show that for certain ``antiresonance'' values of the resonator mass $\mu$, the MiM lattice with Hertzian bead potential (i.e., $\V = \V_{\Hertz}$ in \eqref{original lattice equations}) retains solitary traveling waves.
In a forthcoming paper, Faver, Goodman, and Wright \cite{faver-goodman-wright} do the same for MiM lattices with a potential similar to the Friesecke-Pego class of \ref{fp case}.
We discuss these MiM solitary waves in much greater detail in Section \ref{solitary wave antiresonance section}.
\end{enumerate}

We have chosen the cubic FPUT potential $\V_{\cubic}$ from \eqref{cubic potential} to govern the interaction of the beads in the MiM lattice studied in this paper for two reasons.
First, the simple ``linear + quadratic'' structure\footnote{We expect could add higher-order terms to $\V_{\cubic}$ and use the composition operator techniques from \cite{faver-spring-dimer} to solve the technically more complicated problem that then arises.  For simplicity, however, we omit these terms here. An interesting generalization for future study would be to add nonlinear terms to the bead-resonator force; see Section \ref{future}.} that the spring forces consequently possess puts us in concert with a host of previous papers studying traveling waves in lattices, namely \cite{friesecke-pego1, faver-wright, hoffman-wright, faver-spring-dimer, faver-hupkes}.
Second, we view the cubic FPUT potential as an approximation to the Hertzian potential with precompression.
Taylor-expanding $\V_{\Hertz}'$ from \eqref{Hertzian potential} to second order, one has
\[
\V_{\Hertz}'(r)
\approx \ep_0^p + ar + br^2
\]
for constants $a$ and $b$, and clearly the constant term $\ep_0^p$ is annihilated by the differences in the first equation in both \eqref{original lattice equations} and \eqref{rel disp equations}.
With the bead potential now fixed, we next discuss the selection of the precise monatomic solution from which we will perturb when $\mu$ is small.

\subsection{The monatomic profile}
Rather than restrict ourselves to a Friesecke-Pego \cite{friesecke-pego1} or Friesecke-Wattis \cite{friesecke-wattis} profile, both of which exist for the cubic potential $\V_{\cubic}$, we instead hypothesize the existence of a monatomic solution with some general properties.
First, we need to specify the function spaces in which we work.

For $q > 0$, let 
\[
H_q^r
= \set{f \in H^r}{\cosh^q(\cdot)f \in H^r},
\qquad
\norm{f}_{r,q}
= \norm{f}_{H_q^r}
:= \norm{\cosh^q(\cdot)f}_{H^r}
\]
and set
\[
E_q^r 
:= H_q^r \cap \{\text{even functions}\},
\]
\[
O_q^r
:= H_q^r \cap \{\text{odd functions}\},
\]
and
\[
E_{q,0}^r
:= \set{f \in E_q^r}{\hat{f}(0) = 0},
\]
where $\hat{f}$ is the Fourier transform of $f$ (see Appendix \ref{fourier appendix} for our notations and conventions for Fourier analysis).
These exponentially localized spaces are all Banach spaces under the norm $\norm{\cdot}_{r,q}$; see Appendix C of \cite{faver-dissertation} for some elaborations on their properties.

\begin{hypothesis}\label{hypothesis}
There exists a real number $c$ with $|c| > 1$ such that the following properties hold.

\begin{enumerate}[label={\bf(\roman*)}]

\item
There exists a positive solution $\sigma_c \in E_{q_c}^2$ to the monatomic FPUT traveling wave problem \eqref{monatomic problem} with $\V = \V_{\cubic}$.

\item
The linearization of the monatomic traveling wave problem \eqref{monatomic problem} at $\sigma_c$, which is the operator
\begin{equation}\label{H defn}
\H_cf 
:= c^2f'' - A(1+2\varsigma_c)f,
\end{equation}
is invertible from $E_{q_c}^2 \to E_{q_c,0}^0$.
\end{enumerate}
\end{hypothesis}

When $|c| \gtrsim 1$ and $\sigma_c$ is the Friesecke-Pego \cite{friesecke-pego1} solitary wave, Hoffman and Wright have shown the invertibility of $\H_c$ in Proposition 3.1 of \cite{hoffman-wright}.
Herrmann and Matthies discuss the invertibility of the analogous operator for their singular potential in the high-energy limit in \cite{herrmann-matthies-uniqueness}.
So, we feel there is solid precedent for requiring the invertibility of $\H_c$ in our hypothesis.
The condition that $|c| > 1$ will be important at a technical level deep in the appendices (see the proof of Proposition \ref{eigenprops}, for example), but otherwise it plays no role; we hope that making the existence of a monatomic FPUT solution a hypothesis will put us in further dialogue with existing and future monatomic results, instead of just restricting, say, to the near-sonic case $|c| \gg 1$.

A straightforward bootstrapping argument shows that $\sigma_c \in \cap_{r=0}^{\infty} E_{q_c}^r$.
Less trivial is the following extension of the invertibility of $\H_c$ over a range of weights.
This is stated and proved as Proposition 5.1 in \cite{faver-hupkes}.

\begin{proposition}\label{H prop}
There exist $q_1(c) < q_2(c) \le q_c$ such that for $q \in [q_1(c),q_2(c)]$ and $r \ge 0$, the operator $\H_c$ is invertible from $E_q^{r+2}$ to $E_q^r$.
\end{proposition}

\subsection{Nanopterons and the main result}

Our instinctive intention is to look for traveling wave profiles $(p_1,p_2)$ for the MiM lattice with cubic potential that are small, exponentially localized perturbations of $(\sigma_c,0)$.
However, experience with traveling waves in heterogeneous FPUT lattices, in which the FPUT lattice no longer has identical springs and masses but either alternating masses (the mass dimer or diatomic lattice) or alternating spring potentials (the spring dimer), suggests that material heterogeneity in lattices induces the classical solitary wave to perturb into a {\it{nonlocal}} solitary wave.
Specifically, we refer to the long wave limit studied by Faver and Wright \cite{faver-wright} for the mass dimer and again by Faver \cite{faver-spring-dimer} for the spring dimer, the small mass limit in mass dimers of Hoffman and Wright \cite{hoffman-wright}, and the equal mass limit of Faver and Hupkes \cite{faver-hupkes}, again for the mass dimer.
We also note a variety of numerical studies of nanopeterons in mass and spring dimer lattices, including \cite{gmwz, lustri-porter, vswp, sv, co-ops, lustri}.

The traveling waves for these dimers under the relevant physical limits do not vanish at infinity but instead asymptote to small amplitude periodic ``ripples,'' which themselves are exact traveling waves, too.
In the long wave and small mass problems, a singular perturbation in the traveling wave equations causes the amplitude of the ripple to be small beyond all orders of the relevant small parameter, and thus these nonlocal solitary waves are {\it{nanopterons}}, in Boyd's terminology \cite{boyd}; in the equal mass limit there is no singular perturbation, and the ripples are only, ostensibly, algebraically small.

Furthermore, there is a considerable body of numerical evidence specifically for MiM lattices indicating the existence of nanopterons therein; see the numerical experiments and results discussed in \cite{kev-vain-et-al, woodpile, ksx, xu, vorotnikov-et-al}.
The article \cite{kev-vain-et-al}, for example, suggests that the resonator can act as a ``defect'' to trap some of an MiM traveling wave's energy in localized oscillations. 

\begin{figure}
\[
%
%
\begin{tikzpicture}

\draw[thick,<->] (-7.5,0)--(7.5,0)node[right]{$x$};
\draw[thick,<->] (0,-.5)--(0,2.5);

\def\a{20};
\def\b{1000};
\def\c{20};
\def\d{1.75};
\def\e{75};
\def\f{1.25};

\draw[densely dotted,ultra thick] plot[domain=-7.5:7.5,smooth,samples=500] (\x,{1.8/(1+3*(\x)^2)});

\draw[ultra thick,blue] plot[domain=-7.5:7.5,smooth,samples=500] (\x,
{
\d*exp(-abs(\x)/\e)/(1+(\f*\x)^2)
});

\draw[thick,->] (.25,.45)node[below=-2pt]{$\sigma_c$}--(.5,.8);
\draw[thick,<-] (-.9,.9)--(-1.4,1.4)node[left, inner sep = 0pt]{the nanopteron profile $p_1$};

\draw[densely dotted,line width = .75pt] (6.25,{2/(exp(6.25)+exp(-6.25))}) circle(.5);
\draw[densely dotted, line width = .75pt] (5.75,{2/(exp(5.75)+exp(-5.75))}) -- (5.0253,2.25);
\draw[densely dotted, line width = .75pt] (6.75,{2/(exp(6.75)+exp(-6.75))}) -- (7.4747,2.25);
\draw[densely dotted, line width = .75pt] (6.25,2.5) circle(1.25);

\begin{scope}
\clip (6.25,2.5) circle(1.25);
\node at (6,2.5){\begin{tikzpicture}
\draw[ultra thick,blue] plot[domain=-3:3,samples=100,smooth] (\x,{-.3*cos(10*\x r)});
\end{tikzpicture}};
\end{scope}

\draw [line width = .75pt,decorate,decoration={brace,amplitude=2pt}] (4.75,2.19)--node[midway,left]{amplitude $\sim\O(\mu^{\infty})$}(4.75,2.81);
\draw[line width = .75pt, densely dotted] (4.75,2.19)--(5.36,2.19);
\draw[line width = .75pt, densely dotted] (4.75,2.81)--(5.05,2.81);

\end{tikzpicture}
%
%
%
%
%
%
%
%
%
%
%
\]
\caption{The nanopteron traveling wave profile $p_1$. 
By $\O(\mu^{\infty})$ we mean ``small beyond all orders of $\mu$.''}
\label{fig-nanopteron}
\end{figure}
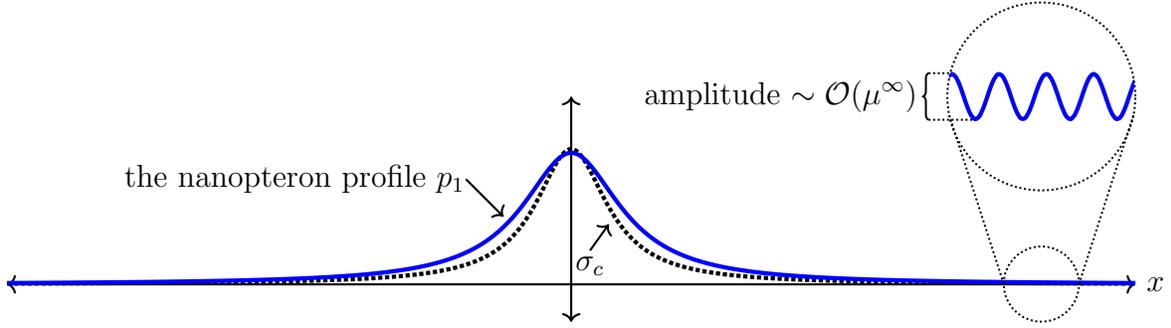

Motivated, then, by this numerical evidence --- in particular the observations that Kevrekidis, Stefanov, and Xu \cite{ksx} made of nanopteron behavior for resonator mass values $\mu$ away from the antiresonance cases --- and the success of Hoffman and Wright in studying the diatomic FPUT small mass problem, we will construct nanopteron traveling waves for the MiM lattice whose bead interactions have the cubic FPUT potential.
Here, then, is our main theorem, proved precisely as Theorem \ref{main nanopteron theorem}.
We sketch the nanopteron profile $p_1$ from this theorem in Figure \ref{fig-nanopteron}.

\begin{theorem}\label{main theorem intro}
Suppose that $|c| > 1$ satisfies Hypothesis \ref{hypothesis} and $\kappa > 0$.
Then there exist $\mu_{\star} >0$ and an open set $\Mu_{c,\kappa} \subseteq (0,\mu_{\star})$ with the following properties.

\begin{enumerate}[label={\bf(\roman*)}]

\item
$0 \in \overline{\Mu_{c,\kappa}}$.

\item
If $\mu \in \Mu_{c,\kappa}$, then there are functions $\varphi_1^{\mu}$, $\varphi_2^{\mu}$, $\Upsilon_1^{\mu}$, $\Upsilon_2^{\mu} \in \Cal^{\infty}(\R)$ such that putting
\begin{equation}\label{intro thm profs}
p_1 := \sigma_c\left(x+\frac{1}{2}\right) +\varphi_1^{\mu}\left(x+\frac{1}{2}\right) + \Upsilon_1^{\mu}\left(x+\frac{1}{2}\right)
\quadword{and}
p_2 := \varphi_2^{\mu}(x) + \Upsilon_2^{\mu}(x)
\end{equation}
gives traveling wave profiles for the MiM equations \eqref{rel disp equations} with cubic FPUT potential \eqref{cubic potential}. 

\item
The functions $\varphi_1^{\mu}$ and $\varphi_2^{\mu}$ are periodic with frequency $\O(\mu^{-1/2})$ and amplitude small beyond all orders of $\mu$; the functions $\Upsilon_1^{\mu}$ and $\Upsilon_2^{\mu}$ are exponentially localized with amplitude $\O(\mu)$.
\end{enumerate}
\end{theorem}

The additional shift by $1/2$ in $p_1$ in \eqref{intro thm profs} is a consequence of a symmetrization choice that we make in our traveling wave problem in \eqref{rhob defn}.
The restriction of $\mu$ to the set $\Mu_{c,\kappa}$ is nontrivial and technical and is discussed at length in Sections \ref{eta2 a eqn motiv} and \ref{solitary wave antiresonance section}.
Roughly, $\Mu_{c,\kappa}$ excludes from the interval $(0,1)$ the antiresonance values of $\mu$ at which \cite{ksx} and \cite{faver-goodman-wright} find solitary waves.

\subsection{Remarks on notation}
We outline some conventions that will be frequently used in the following arguments.

\begin{enumerate}[label=$\bullet$]

\item
We regard $c$ and $\kappa$ as fixed throughout this problem.
All constants and estimates are allowed to depend on both $c$ and $\kappa$.
We occasionally indicate dependence of an operator on $c$ to highlight the relationship to wave speed.

\item
All constants are independent of $\mu$, and if a constant does depend on some other parameter (which is not $c$ or $\kappa$), we indicate that with function notation, e.g., if the constant $C$ depends on the parameters $q$ and $r$, then we write $C = C(q,r)$.

\item
If $x_{\mu}$ is a real number that depends on $\mu$, and if $p > 0$, then we write $x = \O(\mu^p)$ if there is a constant $C(p) > 0$ such that $|x_{\mu}| \le C(p)\mu^p$ for $\mu$ close to 0.
If $x_{\mu} = \O(\mu^p)$ for every $p > 0$, then we write $x_{\mu} = \O(\mu^{\infty})$.

\item
If $\X$ and $\Y$ are normed spaces, we write $\b(\X,\Y)$ for the space of bounded linear operators from $\X$ to $\Y$, and we define $\b(\X) := \b(\X,\X)$.
\end{enumerate}
\section{The Traveling Wave Problem}\label{traveling wave problem section}

\subsection{The traveling wave equations}\label{tw eqns section}
With $\pb = (p_1,p_2)$ and $I_{\mu} := \diag(1,\mu)$, we use the traveling wave ansatz \eqref{tw ansatz} to convert our relative displacement equations \eqref{rel disp equations} to
\begin{equation}\label{pre-symmetry matrix-vector tw}
c^2I_{\mu}\pb'' + \begin{bmatrix*}
-(S^1-2+S^{-1}) &\kappa(S^1-1) \\
\mu(S^{-1}-1) &\kappa(1+\mu)
\end{bmatrix*}\pb
+ \begin{pmatrix*}
-(S^1-2+S^{-1})p_1^2 \\
\mu(S^{-1}-1)p_1^2
\end{pmatrix*}
= 0.
\end{equation}
Motivated by the traveling wave ansatz in \cite{herrmann-matthies-asymptotic}, we can ``symmetrize'' this equation by changing variables to
\begin{equation}\label{rhob defn}
\rhob := \begin{bmatrix*}
S^{-1/2} &0 \\
0 &1
\end{bmatrix*}\pb.
\end{equation}
The extra shift by $1/2$ in our final form of the profiles in \eqref{intro thm profs} is an artifact of this extra change of variables.

Then \eqref{pre-symmetry matrix-vector tw} is equivalent to
\begin{equation}\label{main equation}
\bunderbrace{c^2I_{\mu}\rhob''
+ \D_{\mu}\rhob
+ \D_{\mu}\nl(\rhob,\rhob)}{\G_c(\rhob,\mu)} = 0,
\end{equation}
where
\begin{equation}\label{D-mu defn}
\D_{\mu} := 
\begin{bmatrix*}
-A &\kappa\delta\\
-\mu\delta &\kappa(1+\mu)
\end{bmatrix*},
\qquad
A := S^1-2+S^{-1},
\quadword{and}
\delta := S^{1/2}-S^{-1/2},
\end{equation}
and
\begin{equation}\label{Q defn}
\nl(\rhob,\grave{\rhob})
:= I_0(\rhob.\grave{\rhob})
= \begin{pmatrix*}
\rho_1\grave{\rho_1} \\
0
\end{pmatrix*}.
\end{equation}

\subsubsection{Symmetries}
We observe several useful properties of the mapping $\G_c$.
If $f$ is even and $g$ is odd, then $Af$ and $\delta{g}$ are even, while $Ag$ and $\delta{f}$ are odd.
It follows that if $\rho_1$ is even and $\rho_2$ is odd, then $\G_c(\rhob,\mu)\cdot\e_1$ is even and $\G_c(\rhob,\mu)\cdot\e_2$ is odd, where $\e_1 = (1,0)$ and $\e_2 = (0,1)$.
Next, $Af$ and $\delta{f}$ are ``mean-zero'' in the sense that if $f$ has a Fourier transform $\hat{f}$, then $\hat{Af}(0) = \hat{\delta{f}}(0) = 0$. 
Consequently, $\hat{\G_c(\rhob,\mu)\cdot\e_1}(0) = 0$, as well.
We conclude that $\G_c$ maps
\begin{multline*}
\{\text{even functions}\} \times \{\text{odd functions}\} \times \R \\
\to
\{\text{even mean-zero functions}\} \times \{\text{odd functions}\}.
\end{multline*}

\subsubsection{The $\mu=0$ limit}
We are interested solutions to the traveling wave problem \eqref{main equation} in the limit as $\mu \to 0^+$, so we take $\mu=0$ and find that \eqref{main equation} reduces to the system
\[
\begin{cases}
c^2\rho_1'' - A(\rho_1+\rho_1^2) + \kappa\delta\rho_2 = 0 \\
\kappa\rho_2 = 0.
\end{cases}
\]
The second equation clearly has the solution $\rho_2 = 0$, at which point the first equation becomes
\begin{equation}\label{monatomic tw eqn}
c^2\rho_1'' - A(\rho_1+\rho_1^2) = 0.
\end{equation}
This is the advance-delay differential equation that traveling wave profiles for the monatomic FPUT lattice must satisfy, and by Hypothesis \ref{hypothesis} it has the solution $\rho=\sigma_c$.

\subsection{Linearizing at the monatomic solution}
To perturb from $\sigmab_c$, we make the ansatz 
\begin{equation}\label{original perturbation}
\rhob = \sigmab_c + \varrhob
\end{equation}
for some small, exponentially localized $\varrhob = (\varrho_1,\varrho_2)$.
We then rearrange the equation $\G_c(\sigmab_c + \varrhob,\mu) = 0$ into the system
\begin{equation}\label{naive perturb}
\begin{cases}
\H_c\varrho_1 = \rhs_1(\rhob,\mu) \\
\L_{\mu}\varrho_2 = \rhs_2(\rhob,\mu).
\end{cases}
\end{equation}

The operator $\H_c$ was defined in \eqref{H defn}; the new operator $\L_{\mu}$ is given by
\begin{equation}\label{L-mu defn}
\L_{\mu}f
:= c^2\mu{f}'' +\kappa(1+\mu)f.
\end{equation}
The terms on the right side, absorbed into $\rhs_{c,1}^{\mu}$ and $\rhs_{c,2}^{\mu}$, are all ``small'' in that up to some $\O(1)$ linear operators, they are $\mu$, $\mu\varrhob$, or $\varrhob^{.2} = (\varrho_1,\varrho_2)$.
The precise forms of these operators $\rhs_{c,1}^{\mu}$ and $\rhs_{c,2}^{\mu}$ are not important, as we will not work with exactly these terms in the future.
Hypothesis \ref{hypothesis} provides the invertibility of $\H_c$, so we obtain a fixed point equation for $\varrho_1$ easily:
\[
\varrho_1
= \H_c^{-1}\rhs_1(\rhob,\mu).
\]

However, the operator $\L_{\mu}$ is not invertible.
It is a singularly perturbed, constant coefficient second-order linear differential operator, and so we know that 
\begin{equation}\label{L-mu vanishes on periodics}
\L_{\mu}\sin(\Omega_{\mu}\cdot)
= \L_{\mu}\cos(\Omega_{\mu}\cdot)
= 0,
\end{equation}
where
\begin{equation}\label{Omega-mu defn}
\Omega_{\mu}
:= \sqrt{\frac{\kappa(1+\mu)}{c^2\mu}}.
\end{equation}
We can also view $\L_{\mu}$ as a Fourier multiplier in the sense that 
\[
\hat{\L_{\mu}f}(k)
= (-c^2k^2+\kappa(1+\mu))\hat{f}(k), \ f \in H_q^2.
\]
Importantly, the symbol of $\L_{\mu}$ vanishes at $k=\pm\Omega_{\mu}$, and so
\begin{equation}\label{solv1}
\hat{\L_{\mu}f}(\pm\Omega_{\mu}) = 0.
\end{equation}
Consequently, $\L_{\mu}$ cannot be surjective as an operator $O_q^{r+2} \to O_q^r$, and so we have no hope of inverting $\L_{\mu}$ and converting the second equation in \eqref{naive perturb} into a tractable fixed point problem for $\varrho_2$.

However, while the solvability condition \eqref{solv1} impedes our progress with a classical perturbation attempt to solve the traveling wave problem, it will indicate a different method of resolution.
First, the next proposition shows that $\L_{\mu}$ is injective and provides an exact coercive estimate for its inverse; moreover, we can use \eqref{solv1} to characterize the range of $\L_{\mu}$.
The proof is in Appendix \ref{L prop proof appendix}.

\begin{proposition}\label{L prop}
Let $r \ge 0$ and $\mu$, $q >0$ and set
\[
\Dfrak_{\mu,q}^r 
:= \set{f \in O_q^r}{\hat{f}(\Omega_{\mu}) = 0}.
\]
Then $\L_{\mu}$ is invertible from $O_q^{r+2}$ to $\Dfrak_{\mu,q}^r$.
Denoting the inverse operator from $\Dfrak_{\mu,q}^r$ to $O_q^{r+2}$ by $\L_{\mu}^{-1}$, we have
\begin{equation}\label{L-mu coercive}
\norm{\L_{\mu}^{-1}}_{\b(\Dfrak_{\mu,q}^r,\O_q^{r+j})}
\le C(q)\mu^{-(1+j)/2}, \ j = 0,1,2,
\end{equation}
where $C(q)$ is independent of $\mu$ and $r$, but $C(q) \to \infty$ as $q \to 0^+$.
\end{proposition}

The pointwise evaluation of an exponentially localized function's Fourier transform at the critical value $\Omega_{\mu}$ will be an essential linear functional for our future study. 
We name it and estimate it in the next proposition, whose proof is the same as that of Lemma A.5 in \cite{faver-wright}, which in turn is a modification of the classical Riemann-Lebesgue estimate for the Fourier transform.

\begin{proposition}\label{iota-mu prop}
For $g \in O_q^r$, let
\begin{equation}\label{iota-mu defn}
\iota_{\mu}[g]
:= \int_{-\infty}^{\infty} g(x)\sin(\Omega_{\mu}x) \dx.
\end{equation}
Then for all $q$, $r > 0$, there is $C(q,r) > 0$ such that if $f \in O_q^r$, then
\begin{equation}\label{iota-mu est}
|\iota_{\mu}[f]|
\le C(q,r)\mu^{r/2}.
\end{equation}
\end{proposition}

Using the language of the preceding two propositions, we can then characterize the range of $\L_{\mu}$ via
\begin{equation}\label{solv2}
\L_{\mu}f = g, \ f \in O_q^{r+2}, \ g \in O_q^r
\iff \iota_{\mu}[g] = 0.
\end{equation}
In the next section, we discuss how this solvability condition motivates a revised traveling wave ansatz, the nanopteron.

\subsection{The nanopteron strategy}
The traveling wave problem \eqref{naive perturb} is overdetermined: we must solve the two equations $\varrho_1 = \H_c^{-1}\rhs_1(\varrhob,\mu)$ and $\L_c\varrho_2 = \rhs_2(\varrhob,\mu)$, while the solvability condition \eqref{solv2} demands that the equation $\iota_{\mu}[\rhs_2(\varrhob,\mu)] = 0$ be met as well.
A very similar situation arose in the long wave problem for diatomic FPUT lattices \cite{faver-wright}.
There, the traveling wave problem contained two unknowns, like our $\varrho_1$ and $\varrho_2$, but a third equation, like \eqref{solv2}, had to be met.  
And, like our \eqref{solv2}, the third equation in the long wave problem involved a functional of the form ``integrate against a sinusoid.''

Similar solvability conditions arise in the spring dimer long wave problem \cite{faver-spring-dimer}, and a long wave problem for the capillary-gravity Whitham equation, studied by Johnson and Wright \cite{johnson-wright}.  
Like the prototypical model problem of Amick and Toland \cite{amick-toland}, all of these solvability conditions hinge on the vanishing of an ``integrate against a sinusoid'' functional.
Substantially more complicated solvability conditions appear in the small mass FPUT problem \cite{hoffman-wright} and in the equal mass FPUT problem \cite{faver-hupkes}.
Both still require a functional of the form ``integrate against a function $\gamma$'' to vanish, but the precise function $\gamma$ is no longer an exact sinusoid.
For the small mass FPUT limit, $\gamma$ is an asymptotically sinusoidal ``Jost solution'' of a small perturbation of a Schrodinger operator by nonlocal terms, while in the equal mass limit $\gamma$ is an asymptotically sinusoidal kernel element of an advance-delay differential operator that has exponentially localized, variable coefficient terms.
In the small mass and equal mass FPUT problems, the chief technical challenge (and accomplishment) is, arguably, the construction of the function $\gamma$ that gives the precise solvability condition.

But despite the differences in the exact form of the solvability condition, all of these problems employ a related technique, originally due to Beale's work on capillary-gravity water waves \cite{beale2} and subsequently refined by Amick and Toland for a simpler model problem.
Here is the chief idea.
The linearization of our putative fixed point problem
\begin{equation}\label{putative fp}
\begin{cases}
\varrho_1 = \H_c^{-1}\rhs_1(\varrhob,\mu) \\
\L_c\varrho_2 = \rhs_2(\varrhob,\mu)
\end{cases}
\end{equation}
is the diagonal operator $\diag(\ind,\L_{\mu})$, and, by \eqref{L-mu vanishes on periodics}, we see that this operator annihilates the sinusoidal functions $(0,\sin(\Omega_{\mu}\cdot))$ and $(0,\cos(\Omega_{\mu}\cdot))$, where $\Omega_{\mu}$ is defined in \eqref{Omega-mu defn}.
The operator $\diag(\ind,\L_{\mu})$ is injective from $E_q^{r+2} \times O_q^{r+2}$ to $E_q^r \times O_q^r$ but not surjective due to \eqref{solv2}.

Beale encountered a similar obstacle and first observed that these exact sinusoidal solutions to the linearization extend to exact periodic solutions of the full nonlinear problem \eqref{putative fp}.
The same is true for our small mass MiM problem.
That is, there are solutions to \eqref{putative fp} of the form $\rhob = a\phib_{\mu}^a$ for $\mu$ and $a$ small, where $\phib_{\mu}^a$ is periodic and $a$ is, roughly, an ``amplitude parameter.''
We prove this in Theorem \ref{periodic solutions thm}.
Beale then modified his perturbation ansatz to incorporate these periodic terms, and we do the same, replacing \eqref{original perturbation} with
\[
\rhob 
= \sigmab_c + a\phib_{\mu}^a + \etab,
\]
where $\etab = (\eta_1,\eta_2) \in E_q^r \times O_q^r$ for appropriate $q$ and $r$.
We then have three unknowns --- $a$, $\eta_1$, and $\eta_2$ --- and we parley the two equations in \eqref{putative fp} and the solvability condition into three fixed point equations for these three unknowns.
This is the content of most of Section \ref{nanopteron problem section}.
The construction of an equation for $\eta_1$ is easy, as it is effectively given already in fixed point form as the first equation in \eqref{putative fp}.
In Proposition \ref{A B prop} we spell out exactly how one exploits the solvability condition in conjunction with the second equation in \eqref{putative fp} to summon the separate equations for $\eta_2$ and $a$.
This is all in line with the spirit of nanopteron constructions past in \cite{amick-toland, faver-wright, faver-spring-dimer, johnson-wright, faver-hupkes}.

\subsection{A refined leading order limit}\label{refined leading order limit section}
Before proceeding, we need to change slightly the monatomic solution from which we perturb.
An easy calculation shows
\[
\G_c(\sigmab_c,\mu)
= \begin{pmatrix*}
0 \\
-\mu\delta(\sigma_c+\sigma_c^2)
\end{pmatrix*},
\]
and so $\sigmab_c$ solves $\G_c(\rhob,\mu) = 0$ to $\O(\mu)$.
It turns out that we can, and in fact must, do better.
We will find a function $\varsigmab_c^{\mu}$ that is an $\O(\mu)$ perturbation of $\sigmab_c$ and that solves $\G_c(\rhob,\mu) = 0$ to $\O(\mu^2)$.
We will then replace $\sigmab_c$ in our perturbation ansatz by $\varsigmab_c^{\mu}$.
The need to make this ``refined leading order limit,'' which was also necessary in the small mass FPUT problem, arises for very technical reasons deep in the appendices; see Remark \ref{refined leading order limit remark} for a discussion much later in the right context of why this change matters.

Here is our result, which we prove using the implicit function theorem, exactly as Hoffman and Wright establish Lemma 4.1 in their paper.

\begin{lemma}\label{refined leading order limit thm}
There exists $\mu_{\varsigma} > 0$ such that for all $\mu \in (0,\mu_{\varsigma})$, there is $\zetab_c^{\mu} = (\zeta_{c,1}^{\mu},\zeta_{c,2}^{\mu}) \in \cap_{r=0}^{\infty} E_{q_c}^r \times O_{q_c}^r$ such that if 
\[
\varsigmab_c^{\mu} 
= \begin{pmatrix*} \varsigma_{c,1}^{\mu} \\ \varsigma_{c,2}^{\mu} \end{pmatrix*}
:= \sigmab_c + \mu\zetab_c^{\mu},
\]
then $\G_c(\varsigmab_c^{\mu},\mu) \cdot \e_1 = 0$ and $\G_c(\varsigmab_c^{\mu},\mu) \cdot \e_2 = \O(\mu^2)$.
Moreover, for each $r \ge 0$,
\[
\sup_{\mu \in (0,\mu_{\varsigma})} \norm{\zetab_c^{\mu}}_{r,q_c}
\le C(r).
\]
\end{lemma}
\section{Periodic Solutions}\label{periodic solutions section}

In this section we state a theorem on the existence of periodic solutions to our traveling wave problem \eqref{main equation}.
First, we rescale our traveling wave profiles $\rhob$ to
\begin{equation}\label{periodic ansatz}
\rhob(x) = I^{\mu}\phib(\omega{x}),
\end{equation}
where $I^{\mu} := \diag(\mu,1)$, $\phib(X) = (\phi_1(X),\phi_2(X))$ is $2\pi$-periodic, and $\omega \in \R$.
We assume $\phib \in E_{\per,0}^2 \times O_{\per}^2$, where
\[
H_{\per}^r = \set{f \in L_{\per}^2}{\norm{f}_r < \infty},
\qquad
\ip{f}{g}_r := \sum_{k=-\infty}^{\infty} (1+k^2)^r\hat{f}(k)\overline{\hat{g}}(k),
\qquad
\norm{f}_r := \sqrt{\ip{f}{f}_r},
\]
\[
E_{\per}^r = H_{\per}^r \cap \{\text{even functions}\},
\]
\[
E_{\per,0}^r = \set{f \in E_{\per}^r}{\hat{f}(0) = 0},
\]
and
\[
O_{\per}^r = H_{\per}^r \cap \{\text{odd functions}\}.
\]

Inserting the periodic ansatz \eqref{periodic ansatz} into \eqref{main equation}, we find $\G_c(I^{\mu}\phib(\omega\cdot),\mu) = 0$ if and only if
\begin{equation}\label{periodic problem}
\bunderbrace{c^2\mu\omega^2\phib'' + \D_{\mu}[\omega]I^{\mu}\phib+\D_{\mu}[\omega]\nl(I^{\mu}\phib,I^{\mu}\phib)}{\Phib_{\mu}(\phib,\omega)} = 0.
\end{equation}
Here we define $\D_{\mu}[\omega]$ by replacing any shift operator $S^d$ in the definition of $\D_{\mu}$ in \eqref{D-mu defn} with $S^{d\omega}$.
We will construct solutions to this periodic problem \eqref{periodic problem} using a modification of the ``bifurcation from a simple eigenvalue'' argument of Crandall and Rabinowitz \cite{crandall-rabinowitz} and Zeidler \cite{zeidler} that has already been deployed successfully in the FPUT long wave (Theorem 4.1 in \cite{faver-wright} and Theorem 3.1 in \cite{faver-spring-dimer}), small mass (Theorem 5.1 in \cite{hoffman-wright}), and equal mass (Proposition 3.3 in \cite{faver-hupkes}) problems.
The analysis takes place chiefly at the linear level.

First we examine the kernel of the linearization of \eqref{periodic problem}, and so we study the equation
\[
\bunderbrace{c^2\mu\omega^2\phib'' + \D_{\mu}[\omega]I^{\mu}\phib}{\Gamma_{\mu}[\omega]\phib}
= 0.
\]
Viewed on the Fourier side, a solution $\phib$ to $\Gamma_{\mu}[\omega]\phib = 0$ must satisfy
\begin{equation}\label{ev-ew rel}
c^2\mu\omega^2k^2\hat{\phib}(k) 
= \tD_{\mu}[\omega{k}]I^{\mu}\hat{\phib}(k), \ k \in \Z,
\end{equation}
where $\tD_{\mu}[K] \in \C^{2 \times 2}$ is the matrix
\begin{equation}\label{t-D-mu defn}
\tD_{\mu}[K]
:= \tD_{\mu}[K]
:= \begin{bmatrix*}
-(2\cos(K)-2) &2i\kappa\sin\left(\frac{K}{2}\right) \\
\\
-2i\mu\sin\left(\frac{K}{2}\right) &\kappa(1+\mu)
\end{bmatrix*}.
\end{equation}
If $\hat{\phib}(k) \ne 0$, then \eqref{ev-ew rel} shows that $c^2\mu\omega^2k^2$ must be an eigenvalue of $\tD_{\mu}[\omega{k}]I^{\mu}$ with eigenvector $\hat{\phib}(k)$.
We compute that the eigenvalues of $\tD_{\mu}[K]$ are
\begin{equation}\label{lambda defn}
\lambda_{\mu}^{\pm}(K)
:= \frac{\kappa}{2} + \mu\left(\frac{\kappa}{2} + 1 - \cos(K)\right)
\pm \frac{\sqrt{
P_{\mu}(\cos(K))
}}{2},
\end{equation}
where
\begin{equation}\label{Pmu}
P_{\mu}(t)
:= (\kappa-2t+2)^2\mu^2 + 2\kappa(\kappa+2t-2)\mu + \kappa^2.
\end{equation}

So, if \eqref{ev-ew rel} holds with $\hat{\phib}(k) \ne 0$, we have either
\[
c^2\mu\omega^2k^2 = \lambda_{\mu}^-(\omega{k})
\quadword{or}
c^2\mu\omega^2k^2 = \lambda_{\mu}^+(\omega{k}).
\]
We will show that the first equality never holds, assuming $\phib \in E_{\per,0}^2 \times O_{\per}^2$, while the second holds only when $\omega = \omega_{\mu}/k$, where $k \in \Z\setminus\{0\}$ and $\omega_{\mu}$ is an $\O(\mu^{3/2})$ perturbation of the original critical frequency $\Omega_{\mu}$ from \eqref{Omega-mu defn}.
Specifically, $\omega_{\mu}$ has the following properties.

\begin{lemma}\label{omega-mu lemma}
Let $|c| > 1$.

\begin{enumerate}[label={\bf(\roman*)},ref={(\roman*)}]

\item
There is $\mu_{\omega} > 0$ such that if $0 < \mu < \mu_{\omega}$, there is a unique $\omega_{\mu} > 0$ such that 
\begin{equation}\label{what omega-mu does}
c^2\omega_{\mu}^2-\lambda_{\mu}^+(\omega_{\mu}) = 0.
\end{equation}

\item
$\omega_{\mu} = \O(\mu^{-1/2})$.
More precisely,
\[
\Omega_{\mu}
= \frac{\sqrt{\kappa(1+\mu)}}{c\sqrt{\mu}} 
\le \omega_{\mu}
\le \frac{\sqrt{\lambda_{\mu}^+(\pi)}}{c\sqrt{\mu}}.
\]

\item\label{omega vs Omega}
$\omega_{\mu}-\Omega_{\mu} = \O(\mu^{3/2})$.
\end{enumerate}
\end{lemma}

We use these properties of $\omega_{\mu}$ in a bifurcation argument in Appendix \ref{periodic solutions appendix} to prove the following theorem.

\begin{theorem}\label{periodic solutions thm}
Suppose $|c| > 1$.
There is $\mu_{\per} \in (0,\mu_{\omega})$ and $a_{\per} > 0$ such that for all $0 < \mu < \mu_{\per}$ there are maps
\begin{align*}
[-a_{\per},a_{\per}] &\to \R \colon a \mapsto \omega_{\mu}^a \\
[-a_{\per},a_{\per}] &\to \Cal_{\per}^{\infty} \cap \{\text{even functions}\} \colon a \mapsto \psi_{\mu,1}^a \\
[-a_{\per},a_{\per}] &\to \Cal_{\per}^{\infty} \cap \{\text{odd functions}\} \colon a \mapsto \psi_{\mu,2}^a
\end{align*}
with the following properties.

\begin{enumerate}[label={\bf(\roman*)}]

\item
If
\begin{equation}\label{varphib-mu-a defn}
\phib_{\mu}^a(x)
= \begin{pmatrix*}
\phi_{\mu,1}^a(x) \\
\phi_{\mu,2}^a(x)
\end{pmatrix*}
:= \begin{pmatrix*}
\mu\upsilon_{\mu}\cos(\omega_{\mu}^ax) \\
\sin(\omega_{\mu}^ax)
\end{pmatrix*}
+ \begin{pmatrix*}
\mu\psi_{{\mu},1}^a(\omega_{\mu}^ax) \\
\psi_{{\mu},2}^a(\omega_{\mu}^ax)
\end{pmatrix*},
\end{equation}
where 
\begin{equation}\label{upsilon pseudo-defn}
\upsilon_{\mu} 
= \sin\left(\frac{\omega_{\mu}}{2}\right)\alpha_{\mu}
\end{equation}
and $\alpha_{\mu} = \O(1)$ is defined below in \eqref{alpha-mu defn}, then $\G_c(a\phib_{\mu}^a,\mu) = 0$.

\item
$\omega_{\mu}^0 = \omega_{\mu}$.

\item
$\psi_{\mu,1}^0 = \psi_{\mu,2}^0 = 0$.

\item
For all $r \ge 0$, there is $C_r > 0$ such that if $|a|$, $|\grave{a}| \le a_{\per}$, then
\begin{equation}\label{naked periodic ests}
|\omega_{\mu}^a -\omega_{\mu}^{\grave{a}}|
+ \norm{\psi_{\mu,1}^a - \psi_{\mu,1}^{\grave{a}}}_{W^{r,\infty}}
+ \norm{\psi_{\mu,2}^a-\psi_{\mu,2}^{\grave{a}}}_{W^{r,\infty}}
\le C_r|a-\grave{a}|.
\end{equation}
\end{enumerate}
\end{theorem}
\section{The Nanopteron Problem}\label{nanopteron problem section}

\subsection{Beale's ansatz}\label{beale ansatz section}
We evaluate the traveling wave problem \eqref{main equation} at Beale's ansatz
\begin{equation}\label{Beale ansatz}
\rhob
= \varsigmab_c + a\phib_{\mu}^a + \etab,
\end{equation}
where $|\mu| < \min\{\mu_{\varsigma},\mu_{\per}\}$ and

\begin{enumerate}[label=$\bullet$]

\item 
The function $\varsigmab_c = \sigmab_c + \mu\zetab_c^{\mu}$ is our refined leading order limit from Theorem \ref{refined leading order limit thm}, with $\G_c(\varsigmab_c,\mu)\cdot\e_1 = 0$ and $\G_c(\varsigmab_c,\mu)\cdot\e_2 = \O(\mu^2)$;

\item
The functions $a\phib_{\mu}^a$, $|a| \le a_{\per}$, are the family of periodic solutions from Theorem \ref{periodic solutions thm}, which satisfy $\G_c(a\phib_{\mu}^a,\mu) = 0$;

\item
And $\etab = (\eta_1,\eta_2) \in E_q^r \times O_q^r$, with $q$ and $r$ to be determined later.
\end{enumerate}

After inserting Beale's ansatz \eqref{Beale ansatz} into the traveling wave problem, we use the linearity of $\D_{\mu}$, the bilinearity of $\nl$, and the properties of the terms of the ansatz stated above to obtain the equation 
\begin{equation}\label{nano syst orig}
c^2I_{\mu}\etab'' + \D_{\mu}\etab + 2\D_{\mu}\nl(\varsigmab_c,\etab)
= -\G_c(\varsigmab_c,\mu)
-2a\D_{\mu}\nl(\varsigmab_c,\phib_{\mu}^a)
-2a\D_{\mu}\nl(\phib_{\mu}^a,\etab)
-\D_{\mu}\nl(\etab,\etab).
\end{equation}
Now we expose the components of this system.
We use the definitions of $\D_{\mu}$ and $\nl$ from \eqref{D-mu defn} and \eqref{Q defn}, respectively, and the results that $\G_c(a\phib_{\mu}^a,\mu) = \G_c(\varsigmab_c,\mu) \cdot \e_1 = 0$ to rewrite \eqref{nano syst orig} as
\begin{equation}\label{nano syst comps}
\begin{cases}
\H_c\eta_1
= h_1^{\mu}(\etab,a) + h_2^{\mu}(\etab,a) + h_3^{\mu}(\etab,a) + h_4^{\mu}(\etab,a) - \kappa\delta\eta_2 \\
\\
\L_{\mu}\eta_2
= \ell_1^{\mu}(\etab,a) + \ell_2^{\mu}(\etab,a) + \ell_3^{\mu}(\etab,a) + 2a\mu\delta[\varsigma_{c,1}^{\mu}\phi_{\mu,1}^a],
\end{cases}
\end{equation}
where
\begin{equation}\label{h ell terms}
\begin{aligned}
h_1^{\mu}(\etab,a) &:= 2\mu{A}(\zeta_{c,1}^{\mu}\eta_1)
\quad&\ell_1^{\mu}(\etab,a) &:= -\G_c(\varsigmab_c,\mu)\cdot\e_2 \\
h_2^{\mu}(\etab,a) &:= 2aA(\varsigma_{c,1}^{\mu}\phi_{\mu,1}^a)
\quad&\ell_2^{\mu}(\etab,a) &:= \mu\delta[\eta_1+\varsigma_{c,1}^{\mu}\eta_1+\eta_1^2]\\
h_3^{\mu}(\etab,a) &:= 2aA(\phi_{\mu,1}^a\eta_1)
\quad&\ell_3^{\mu}(\etab,a) &:= 2a\mu\delta(\phi_{\mu,1}^a\eta_1)\\
h_4^{\mu}(\etab,a) &:= A(\eta_1^2).
\end{aligned}
\end{equation}

The operators $\H$ and $\L_{\mu}$ are the same as they were in \eqref{H defn} and \eqref{L-mu defn} from our initial attempt at a perturbation problem.
One should view the seven terms $h_1^{\mu}$, $h_2^{\mu}$, $h_3^{\mu}$, $\ell_1^{\mu},\ldots,\ell_4^{\mu}$ as ``small.''
For example, $h_1^{\mu}$ is roughly $\mu\eta_1$, $h_2^{\mu}$ is roughly $\mu{a}$, when we recall from \eqref{varphib-mu-a defn} that $\phi_{\mu,1}^a$ bears a factor of $\upsilon_{\mu} = \O(\mu)$, and $h_3^{\mu}$ and $h_4^{\mu}$ are quadratic in $a$ and $\mu$.
We choose to leave the last terms in the $\eta_1$ and $\eta_2$ equations explicit for now, as they will require more modifications shortly.

\subsection{The fixed point equations for $\eta_2$ and $a$: motivation}\label{eta2 a eqn motiv}
We motivate our construction of these equations by thinking about how the inverse function theorem could apply to our situation.
The second equation in \eqref{nano syst comps} reads
\[
\L_{\mu}\eta_2
= \bunderbrace{\sum_{k=1}^3 \ell_k^{\mu}(\etab,a) + 2a\mu\delta[\varsigma_{c,1}^{\mu}\phi_{\mu,1}^a]}{\rhs_{\mu}(\etab,a)}.
\]
If this is to hold, the solvability condition \eqref{solv2} requires
\begin{equation}\label{solv3}
\iota_{\mu}[\rhs_{\mu}(\etab,a)]
= 0.
\end{equation}
The inverse function theorem will permit us to solve \eqref{solv3} for $a$ in terms of $\etab$ if 
\begin{equation}\label{inverse function theorem attempt}
\frac{\partial}{\partial{a}}\iota_{\mu}\left[\sum_{k=1}^3 \ell_k^{\mu}(\etab,a)+2a\mu\delta[\varsigma_{c,1}^{\mu}\phi_{\mu,1}^a]\right]\bigg|_{\etab=0,a=0}
= 2\mu\iota_{\mu}\big[\delta[\varsigma_{c,1}^{\mu}\phi_{\mu,1}^0]\big]
\ne 0.
\end{equation}
This calculation holds because the terms $\ell_1^{\mu}$ and $\ell_2^{\mu}$ are constant in $a$, while $\ell_3^{\mu}$ is quadratic in $a$ and $\etab$; this calculation, though, is purely formal because we never established the differentiability of the map $a \mapsto \phib_{\mu}^a$ as part of Theorem \ref{periodic solutions thm}.
And so we will not use the inverse function theorem here, not least because we would require further results on the periodic functions, but also because we would need some highly quantitative modification of that classical theorem to handle the $\mu$ parameter. 

Instead, we modify the second equation in \eqref{nano syst comps} based on what \eqref{inverse function theorem attempt} suggests the leading order term in $a$ should be.
Observe that
\begin{equation}\label{ell-4 ell-5 defns}
2a\mu\delta[\varsigma_{c,1}^{\mu}\phi_{\mu,1}^a]
= 2a\mu\delta[(\sigma_c+\mu\zeta_{c,1}^{\mu})\phi_{\mu,1}^a]
= 2a\mu\delta[\sigma_c\phi_{\mu,1}^0] + \bunderbrace{2a\mu\delta[\sigma_c(\phi_{\mu,1}^a-\phi_{\mu,1}^0)]}{\ell_4^{\mu}(\etab,a)}+\bunderbrace{2a\mu^2\delta(\zeta_{c,1}^{\mu}\phi_{\mu,1}^a)}{\ell_5^{\mu}(\etab,a)}.
\end{equation}
We think of $\ell_4^{\mu}$ as being roughly quadratic in $a$, while $\ell_5^{\mu}$ has the very small factor $\mu^2$.

Theorem \ref{periodic solutions thm} allows us to calculate
\[
\phi_{\mu,1}^0
= \mu\upsilon_{\mu}\cos(\omega_{\mu}\cdot),
\]
where $\upsilon_{\mu} = \O(1)$ is defined in \eqref{upsilon pseudo-defn}.
Recall from Lemma \ref{omega-mu lemma} that $\omega_{\mu}$ is an $\O(\mu^{3/2})$ perturbation of the original critical frequency $\Omega_{\mu}$, defined in \eqref{Omega-mu defn}.
For reasons that will be clearer in the proper context later, specifically Remark \ref{why the frequency change?}, we want to make one further adjustment to permit $\Omega_{\mu}$ to play more of a starring role.
We rewrite
\begin{equation}\label{return to original frequency}
2a\mu\delta[\sigma_c\phi_{\mu,1}^0]
= \bunderbrace{2a\mu^2\upsilon_{\mu}\delta[\sigma_c\cos(\Omega_{\mu}\cdot)]}{-a\chi_{\mu}} 
+ \bunderbrace{2a\mu^2\upsilon_{\mu}\delta\big[\sigma_c\big(\cos(\omega_{\mu}\cdot)-\cos(\Omega_{\mu}\cdot)\big)\big]}{\ell_6^{\mu}(\etab,a)}.
\end{equation}

That is, we have set
\begin{equation}\label{chi-mu defn}
\chi_{\mu}
:= -2\mu^2\upsilon_{\mu}\delta[\sigma_c\cos(\Omega_{\mu}\cdot)].
\end{equation}
Combining \eqref{ell-4 ell-5 defns} and \eqref{return to original frequency}, we see that the second equation in \eqref{nano syst comps} is equivalent to
\begin{equation}\label{eta2 eqn1}
\L_{\mu}\eta_2 + a\chi_{\mu} = \sum_{k=1}^6 \ell_k^{\mu}(\etab,a).
\end{equation}
We pause to isolate some precise estimates on $\chi_{\mu}$.

\begin{lemma}\label{chi-mu rq est lemma}
For each $r \ge 0$ and $q \in [0,q_c]$, we have $\chi_{\mu} \in E_q^r$, and there is $C(q,r) > 0$ such that 
\begin{equation}\label{chi-mu est}
\norm{\chi_{\mu}}_{r,q} 
\le C(q,r)\mu^{2-r/2}.
\end{equation}
\end{lemma}

Now we sketch how we obtain the equation for $a$.
Recall that the functional $\iota_{\mu}$, defined in \eqref{iota-mu defn}, annihilates the range of $\B_{\mu}$.
Applying this functional to \eqref{eta2 eqn1}, we have
\[
a\iota_{\mu}[\chi_{\mu}]
= \sum_{k=1}^6\iota_{\mu}[\ell_k^{\mu}(\etab,a)].
\]
Consequently, we want to write
\begin{equation}\label{a eqn prelim}
a
= \frac{1}{\iota_{\mu}[\chi_{\mu}]}\sum_{k=1}^6\iota_{\mu}[\ell_k^{\mu}(\etab,a)].
\end{equation}
Of course, this is valid only if $\iota_{\mu}[\chi_{\mu}] \ne 0$.

In Appendix \ref{chi-id proof appendix}, we establish the following identity for $\chi_{\mu}$:
\begin{equation}\label{chi-id}
\iota_{\mu}[\chi_{\mu}]
= \mu^2\sin^2\left(\frac{\Omega_{\mu}}{2}\right)\alpha_{\mu}\norm{\sigma_c}_{L^1} + \O(\mu^{7/2}),
\end{equation}
where $\alpha_{\mu} = \O(1)$ is defined precisely in \eqref{alpha-mu defn} and $\norm{\sigma_c}_{L^1} > 0$.
So, the denominator in \eqref{a eqn prelim} is at least $\O(\mu^2)$, which is small.

However, the denominator can be smaller, for it is possible to have $\sin(\Omega_{\mu}/2) = 0$.
Specifically, from the definition of $\Omega_{\mu}$ in \eqref{Omega-mu defn} this happens if and only if 
\begin{equation}\label{mu-n defn}
\mu
= \frac{\kappa}{4\pi^2n^2c^2-\kappa}
=: \mu_n, \ n \in \N.
\end{equation}
One can check that these $\mu_n$ are precisely the ``antiresonance'' mass ratios that Kevrekidis, Stefanov, and Xu encounter in \cite{ksx} and for which they produce solitary traveling waves in their Hertzian MiM lattice; see Section \ref{solitary wave antiresonance section} momentarily for a discussion of how these $\mu_n$ relate to solitary waves in our MiM lattice with cubic FPUT potential.

If we attempt to run our nanopteron program including these $\mu_n$, we will encounter an $\O(\mu^{-7/2})$ factor in the equation \eqref{a eqn prelim} for $a$, and this is simply too large for any reasonable estimate.
Another way to view this complication is that if $\mu=\mu_n$, then the formal derivative in \eqref{inverse function theorem attempt} could be $\O(\mu^{7/2})$, which is to say that it could be 0, in which case the inverse function theorem fails to solve \eqref{solv3}.
So, we exclude these potentially pernicious values $\mu_n$ from consideration by restricting $\mu$ to the set
\begin{equation}\label{Mu-c-kappa defn}
\Mu_{c,\kappa}
:= \set{\mu \in (0,\mu_{\per})}{\left|\sin\left(\frac{\Omega_{\mu}}{2}\right)\right| > \frac{1}{2}}.
\end{equation}
Upon doing so, the formula \eqref{chi-id} implies the next lemma.

\begin{lemma}\label{chi lower bound lemma}
There exist $\mu_{\chi} \in (0,\mu_{\per})$ and $C_{\chi} > 0$ such that if $\mu \in \Mu_{c,\kappa} \cap (0,\mu_{\chi})$, then
\begin{equation}\label{iota-mu on chi-mu est}
C_{\chi}\mu^2 
< |\iota_{\mu}[\chi_{\mu}]|.
\end{equation}
\end{lemma}

\begin{remark}\label{why the frequency change?}
Had we not added zero and made the expansion \eqref{return to original frequency}, all the subsequent calculations above would still go through, with $\Omega_{\mu}$ replaced by $\omega_{\mu}$ throughout.
We would still need to restrict $\mu$ to a set like $\Mu_{c,\kappa}$, as the continuity and unboundedness of the map $\mu \mapsto \omega_{\mu}$, which one can extract from \eqref{what omega-mu does}, imply that $\sin(\omega_{\mu}/2) = 0$ for a countable number of $\mu$.
However, we do not have an exact formula for $\omega_{\mu}$, only the implicit relation \eqref{what omega-mu does}.  
So, working with the original critical frequency $\Omega_{\mu}$ and its accompanying antiresonance mass ratios $\mu_n$ is much more transparent.
Furthermore, as we discuss in the next section, it is possible to divine the traveling wave behavior of the MiM lattice at these $\mu_n$, too.
\end{remark}

\subsection{Traveling waves at the antiresonance values for $\mu$}\label{solitary wave antiresonance section}
Hoffman and Wright also had to exclude a countable number of mass ratios $\mu$ from their nanopteron program in \cite{hoffman-wright} for exactly the same reasons that led to our restriction to $\Mu_{c,\kappa}$.
Such a restriction is not an egregious omission of small values of $\mu$; since $\kappa/((2n+1)^2\pi^2-\kappa) \in \Mu_{c,\kappa}$, we have $0 \in \overline{\Mu_{c,\kappa}}$.
Moreover, the lower bound of $1/2$ in the definition of $\Mu_{c,\kappa}$ is arbitrary; we could fix any $b \in (0,1)$ and replace $1/2$ with $b$, and still our nanopteron program would run smoothly (with $C_{\chi}$ now depending, arguably badly, on $b$).

All this is to say that, effectively, we can produce nanopteron traveling waves for $\mu$ arbitrarily close to the values $\mu_n$ from \eqref{mu-n defn}.
But what happens when $\mu=\mu_n$?
Our arguments here certainly do not extend to these values.
It turns out that the MiM lattice bears genuine solitary traveling waves when $\mu = \mu_n$, just as the MiM lattice with Hertzian potential does.
This is proved by Faver, Goodman, and Wright in the forthcoming \cite{faver-goodman-wright}.

The crux of their argument is that when $\mu=\mu_n$, it is possible to rewrite the original traveling wave problem \eqref{main equation} as a single equation in the unknown $\rho_1$, to which a more routine perturbation argument can be applied.
This rewriting, in turn, hinges on viewing the second component of the traveling wave problem \eqref{main equation}, which is
\begin{equation}\label{second component for sol}
\L_{\mu_n}\rho_2 
= \mu_n\delta(\rho_1+\rho_1^2),
\end{equation}
on the Fourier side as 
\[
-(c^2\mu_n{k}^2-\kappa(1+\mu_n))\hat{\rho}_2(k)= 2i\mu_n\sin\left(\frac{k}{2}\right)\left(\hat{\rho}_1 + \hat{\rho_1^2}\right).
\]
The function
\[
\tM_n(k) := -\frac{2i\mu_n\sin(k/2)}{c^2\mu_nk^2-\kappa(1+\mu)}
\]
is analytic on $\C$ with removable singularities at $k = \pm\Omega_{\mu_n}$.
Consequently, Fourier multiplier theory due to Beale (Lemma \ref{beale fm} in this paper) implies that the Fourier multiplier $\M_n$ with symbol $\tM_n$ maps $E_q^r$ to $E_q^{r+2}$, and so \eqref{second component for sol} is equivalent to
\[
\rho_2 = \M_n(\rho_1+\rho_1^2).
\]
Then one feeds this equation into the first component of \eqref{main equation} to obtain an equation solely in $\rho_1$. 
It is also possible to show that the operator norm of $\M_n$ is $\O(1/n)$, and thereby hangs the perturbation argument.
We emphasize, though, that the whole procedure very much depends on taking $\mu = \mu_n$ from \eqref{mu-n defn}, so that the zeros of the denominator of $\tM_n$ cancel with those of the numerator.
This, we feel, justifies our insistence on making $\Omega_{\mu}$ explicit in the nanopteron program via \eqref{return to original frequency}.

\subsection{The fixed point equations for $\eta_2$ and $a$: construction}

Now that we have controlled our permitted values of $\mu$, we are ready to state precisely how we extract two fixed point equations, one for $\eta_2$ and one for $a$, from the single equation \eqref{eta2 eqn1}.

\begin{proposition}\label{A B prop}
Let $q > 0$, $r \ge 0$, and $\mu \in \Mu_{c,\kappa} \cap (0,\mu_{\chi})$.
There exist a bounded linear functional $\A_{\mu} \colon O_q^r \to \R$ and a bounded operator $\B_{\mu} \colon O_q^r \to O_q^{r+2}$ such that, given $f \in O_q^{r+2}$, $a \in \R$, and $g \in O_q^r$, we have
\begin{equation}\label{model solv}
\L_{\mu}f + a\chi_{\mu} = g
\end{equation}
if and only if
\[
a = \A_{\mu}g 
\quadword{and}
f = \B_{\mu}g.
\]
Moreover, the operators $\A_{\mu}$ and $\B_{\mu}$ satisfy
\begin{equation}\label{A-mu B-mu ests}
\norm{\A_{\mu}}_{\b(O_q^r,\R)}
\le C(q,r)\mu^{-2+r/2}
\quadword{and}
\norm{\B_{\mu}}_{\b(O_q^r,O_q^{r+j})}
\le C(q,r)\mu^{-(j+1)/2}, \ j = 0,1,2.
\end{equation}
\end{proposition}

\begin{proof}
Suppose that \eqref{model solv} holds and apply the functional $\iota_{\mu}$ from \eqref{iota-mu defn} to both sides of this equation.
Then we can use Lemma \ref{chi lower bound lemma} to solve for $a$:
\begin{equation}\label{a-eqn abstr}
a 
= \frac{\iota_{\mu}[g]}{\iota_{\mu}[\chi_{\mu}]}
=: \A_{\mu}[g].
\end{equation}
We combine the estimates of \eqref{iota-mu est} in Proposition \ref{iota-mu prop} on $\iota_{\mu}$ and \eqref{iota-mu on chi-mu est} in Lemma \ref{chi lower bound lemma} on $\iota_{\mu}[\chi_{\mu}]$ to prove the estimate for $\A_{\mu}$ in \eqref{A-mu B-mu ests}.

Now set
\[
\P_{\mu}g 
:= g-\frac{\iota_{\mu}[g]}{\iota_{\mu}[\chi_{\mu}]}\chi_{\mu}.
\]
The previous estimate on $\iota_{\mu}[\chi_{\mu}]$ and the estimates on $\iota_{\mu}$ from in Proposition \ref{iota-mu prop} and on $\chi_{\mu}$ from Lemma \ref{chi-mu rq est lemma} cancel each other out to imply
\begin{equation}\label{P-mu operator norm}
\norm{\P_{\mu}}_{\b(O_q^r)} \le C(q,r).
\end{equation}
Next, substituting $a = \A_{\mu}[g]$ into \eqref{model solv} implies
\begin{equation}\label{L=P}
\L_{\mu}f = \P_{\mu}g.
\end{equation}

A direct calculation shows $\iota_{\mu}[\P_{\mu}g] = 0$, so, by the solvability condition \eqref{solv2}, $\P_{\mu}g$ is in the range of $\L_{\mu}$ in $O_q^r$.  
Using the notation of Proposition \ref{L prop} for the inverse of $\L_{\mu}$ from its range to $O_q^{r+2}$, the equality \eqref{L=P} implies
\[
f 
= \L_{\mu}^{-1}\P_{\mu}g 
=: \B_{\mu}g.
\]
Then the coercive estimate \eqref{L-mu coercive} on $\L_{\mu}$ in Proposition \ref{L prop} and \eqref{L=P} imply
\[
\norm{\B_{\mu}g}_{r+j,q} 
\le C(r,q)\mu^{-(j+1)/2}\norm{g}_{r,q}.
\]

Conversely, that $f$, $a$, and $g$ satisfy \eqref{model solv} if $a = \A_{\mu}[g]$ and $f = \B_{\mu}g$ follows from a straightforward algebraic rearrangement of the definitions of $\A_{\mu}$ and $\B_{\mu}$.
\end{proof}

We return to our modified equation \eqref{eta2 eqn1} for $\eta_2$ and $a$.
We apply Proposition \ref{A B prop} with $f=\eta_2$ and $g = \medsum_{k=1}^6 \ell_k^{\mu}(\etab,a)$ to see that \eqref{eta2 eqn1} is equivalent to the pair of fixed point equations
\begin{equation}\label{N2 N3 eqn defn}
\eta_2 
= \sum_{k=1}^6 \B_{\mu}\ell_k^{\mu}(\etab,a)
=: \nano_2^{\mu}(\etab,a)
\quadword{and}
a = \sum_{k=1}^6 \A_{\mu}[\ell_k^{\mu}(\etab,a)]
=: \nano_3^{\mu}(\etab,a).
\end{equation}

\subsection{Construction of the fixed point equation for $\eta_1$}
Last, we must construct an equation for $\eta_1$.
The first equation in \eqref{nano syst comps} reads
\begin{equation}\label{eta1 eqn1}
\H_c\eta_1
= \sum_{k=1}^4 h_k^{\mu}(\etab,a) -\kappa\delta\eta_2.
\end{equation}
By Proposition \ref{H prop} we may invert $\H_c$, so that \eqref{eta1 eqn1} is equivalent to
\begin{equation}\label{eta1 eqn2}
\eta_1
= \sum_{k=1}^4 \H_c^{-1}h_k^{\mu}(\etab,a)-\kappa\H_c^{-1}\delta\eta_2.
\end{equation}
Next, we replace\footnote{We could have made an additional change of variables in Section \ref{tw eqns section} that would have converted our symmetrized traveling wave problem \eqref{main equation} into one whose linear part is genuinely diagonal at $\mu=0$; Hoffman and Wright \cite{hoffman-wright} do this in their small mass problem (see Section 1.3 of that paper).
This would result in what we call $h_5^{\mu}$ in \eqref{h5} being an $\O(\mu)$ term from the start, in which case replacing $\eta_2$ with $\nano_2^{\mu}$ would be unnecessary.
However, the new ``almost diagonal'' version of \eqref{main equation} and the resulting nanopteron system after Beale's ansatz would have a number of more complicated and opaque terms, and so we feel it is not worth the effort making this extra change of variables here.} $\eta_2$ with its fixed point equation $\eta_2 = \nano_2^{\mu}(\etab,a)$ and set
\begin{equation}\label{h5}
h_5^{\mu}(\etab,a)
:= -\kappa\delta\nano_2^{\mu}(\etab,a)
\quadword{and}
\nano_1^{\mu}(\etab,a)
:= \sum_{k=1}^5 \H_c^{-1}h_k^{\mu}(\etab,a).
\end{equation}
Then the three equations \eqref{N2 N3 eqn defn} and \eqref{eta1 eqn2} are equivalent to the fixed point problem
\begin{equation}\label{nano fp}
(\etab,a)
= \big(\nano_1^{\mu}(\etab,a),\nano_2^{\mu}(\etab,a),\nano_3^{\mu}(\etab,a)\big)
=: \nanob_{\mu}(\etab,a).
\end{equation}

\subsection{Solution of the nanopteron fixed point problem \eqref{nano fp}}\label{nano soln}

We will use the following fixed point lemma to solve our main problem; it appears as Theorem 4.10 in \cite{johnson-wright}.

\begin{lemma}\label{nano fp lemma}
Let $\X_0$ and $\X_1$ be reflexive Banach spaces with $\X_1 \subseteq \X_0$.
For $r > 0$, let $\Bfrak(r) := \set{x \in \X_1}{\norm{x}_{\X_1} \le r}$.
Suppose that for some $r_0 > 0$, there is a map $\nano \colon \Bfrak(r_0) \to \X_1$ with the following properties.

\begin{enumerate}[label={\bf(\roman*)}]

\item
$
\norm{x}_{\X_1} \le r_0 
\Longrightarrow \norm{\nano(x)}_{\X_1} \le r_0.
$

\item
There exists $\alpha \in (0,1)$ such that 
\[
\norm{x}_{\X_0}, \norm{\grave{x}}_{\X_1} \le r_0
\Longrightarrow
\norm{\nano(x)-\nano(\grave{x})}_{\X_0} 
\le \alpha\norm{x-\grave{x}}_{\X_0}.
\]
\end{enumerate}
Then there exists a unique $x_{\star} \in \X_1$ such that $\norm{x_{\star}}_{\X_1} \le r_0$ and $x_{\star} = \nano(x_{\star})$.
\end{lemma}

First we specify the Banach spaces (really, Hilbert spaces) in which we will work. 
With $q_1(c)$ and $q_2(c)$ from Proposition \ref{H prop}, fix $q_{\star}(c) \in (q_1(c),q_2(c))$ and take $\qbar_{\star}(c) \in (q_1(c),q_{\star}(c))$.
This restriction ensures that $\H_c$ is invertible from $E_{q_{\star}(c)}^{r+2} \to E_{q_{\star}(c),0}^r$ and $E_{\qbar_{\star}(c)}^{r+2} \to E_{\qbar_{\star}(c),0}^r$ for all $r \ge 0$.
Now let
\[
\X_r
:= \begin{cases}
E_{\qbar_{\star}(c)}^3 \times O_{\qbar_{\star}(c)}^3 \times \R, \ r = 0 \\
\\
E_{q_{\star}(c)}^{r+2} \times O_{q_{\star}(c)}^{r+2} \times \R, \ r \ge 1.
\end{cases}
\]

Next, for $r \ge 1$ and $\mu$, $\tau > 0$, set
\begin{equation}\label{Uhood}
\U_{\mu,\tau}^r
:= \set{(\etab,a) \in \X^r}{\norm{\etab}_{r,q_{\star}(c)} \le \tau\mu^{1/2}, \ |a| \le \tau\mu^{r/2-1}}.
\end{equation}
The following lemma, proved in Appendix \ref{proof of main workhorse lemma appendix}, is instrumental in our application of Lemma \ref{nano fp lemma}.

\begin{lemma}\label{main workhorse lemma}
There exist $\mu_{\star} \in (0,\min\{\mu_{\varsigma},\mu_{\chi}\})$ and  $\tau_{\star} > 0$ such that the following hold for all $\mu \in \Mu_{c,\kappa} \cap (0,\mu_{\star})$.

\begin{enumerate}[label={\bf(\roman*)},ref={(\roman*)}]

\item\label{main workhorse map}
If $(\etab,a) \in \U_{\mu,\tau_{\star}}^1$, then $\nanob_{\mu}(\etab,a) \in \U_{\mu,\tau_{\star}}^1$.

\item\label{main workhorse lip}
If $(\etab,a)$, $(\grave{\etab},\grave{a}) \in \U_{\mu,\tau_{\star}}^1$, then 
\[
\norm{\nanob_{\mu}(\etab,a)-\nanob_{\mu}(\grave{\etab},\grave{a})}_{\X_0}
\le \frac{1}{2}\norm{(\etab,a)-(\grave{\etab},\grave{a})}_{\X_0}.
\]

\item\label{main workhorse boot}
For any $\tau > 0$ and $r \ge 1$, if $(\etab,a) \in \U_{\mu,\tau_{\star}}^1 \cap \U_{\mu,\tau}^r$, then there is $\ttau > 0$ such that $(\etab,a) \in \U_{\mu,\ttau}^{r+1}$.
\end{enumerate}
\end{lemma}

Fixing $r=3$, parts \ref{main workhorse map} and \ref{main workhorse lip} of Lemma \ref{main workhorse lemma} allow us to use Lemma \ref{nano fp lemma} to find a unique $(\etab_{\mu},a_{\mu}) \in \U_{\mu,\tau_{\star}}^1$ such that $(\etab_{\mu},a_{\mu}) = \nanob_{\mu}(\etab_{\mu},a_{\mu})$.
Part \ref{main workhorse boot} then shows that $\etab_{\mu}$ is smooth.  
Since $a_{\mu} = \nano_3^{\mu}(\etab_{\mu},a_{\mu})$, we recall the definition of $\nano_3^{\mu}$ in \eqref{N2 N3 eqn defn} and the estimates on the functional $\A_{\mu}$ from Proposition \ref{A B prop} to conclude that $a_{\mu}$ is small beyond all orders of $\mu$.
We summarize these results formally in our main theorem below.

\begin{theorem}\label{main nanopteron theorem}
Let $|c| > 1$, $\kappa > 0$, and $\mu \in \Mu_{c,\kappa} \cap (0,\mu_{\star})$.
Then there exists a unique $(\etab_{\mu},a_{\mu}) \in \U_{\mu,\tau_{\star}}^1$ such that $(\etab_{\mu},a_{\mu}) = \nanob_{\mu}(\etab_{\mu},a_{\mu})$.
Moreover, $\etab_{\mu} \in \cap_{r=0}^{\infty} E_{q_{\star}(c)}^r \times O_{q_{\star}(c)}$ and there are constants $C_r > 0$ such that for each $r \ge 0$, we have
\begin{equation}\label{main nano ests}
\norm{\etab_{\mu}}_{r,q_{\star}(c)}
\le C_r\mu
\quadword{and}
|a_{\mu}| \le C_r\mu^r.
\end{equation}
\end{theorem}

\subsection{Additional questions and considerations}\label{future}

The paper \cite{faver-goodman-wright} also considers the ``stiff internal spring'' limit for MiM lattices.
In this limit, one tries to solve the traveling wave problem for the lattice equations \eqref{original lattice equations} in the limit as $\kappa \to \infty$.
After defining new traveling wave profiles by rescaling those from \eqref{rhob defn} as $\theta_1 = \rho_1$ and $\theta_2 = \rho_2/\kappa$, the traveling wave problem for $\theta_1$ and $\theta_2$ becomes singularly perturbed like \eqref{main equation}.
Faver, Goodman, and Wright show that for a countable number of ``antiresonance'' values of $\kappa$ that diverge to $\infty$, the resulting system has solitary wave solutions.
We expect that the techniques of this paper will carry over to show the existence of nanopteron traveling waves in the stiff internal spring limit for spring constants $\kappa$ bounded away from the antiresonance values.
For that matter, a singularly perturbed stiff spring limit could be posed for the spring dimer FPUT lattice in \cite{faver-spring-dimer}, analogous to the small mass FPUT limit in \cite{hoffman-wright}.

The nanopteron estimates in \eqref{main nano ests} only establish that the amplitude $a_{\mu}$ of the ripple is small beyond all algebraic orders of $\mu$ and do not provide a lower bound on $a_{\mu}$, so that, conceivably, we could have $a_{\mu} = 0$ for some $\mu$, and our nanopteron would really be a solitary wave.
There are a number of numerical studies indicating that this may occur in diatomic FPUT lattices, see \cite{vswp, sv, lustri-porter, lustri}.
Ideally, we would also have an exponential upper bound of the form $|a_{\mu}| \le e^{-C_{\star}/\mu}$ for some $C_{\star} = \O(1)$, which is what Sun and Shen achieve for the capillary-gravity water wave problem \cite{sun-shen}.
The same question can be posed for the amplitude of the ripples in the nanopteron papers \cite{faver-wright, hoffman-wright, faver-spring-dimer, johnson-wright}. 

Additionally, since we now know that solitary waves exist for the antiresonance mass ratios $\mu_n$ from \eqref{mu-n defn}, another natural question is the value of the limit $\tst{\lim_{\mu \to \mu_n} a_{\mu}}$ for any fixed $n$.
Does this limit equal zero?
Of course, for this limit even to make sense, we must recall that we have found nanopterons with ripple amplitude $a_{\mu}$ where $\mu$ was restricted to the set $\Mu_{c,\kappa}$; as discussed at the start of Section \ref{solitary wave antiresonance section}, we could redefine $\Mu_{c,\kappa}$ to allow a more generous lower bound on $|\sin(\Omega_{\mu}/2|$ and thus permit membership of a larger range of $\mu$ closer and closer to the $\mu_n$.

Resolution of any of these ripple questions will hinge on better knowledge of the oscillatory integral operator $\iota_{\mu}$ from \eqref{iota-mu defn}, which is the critical operator in the definition \eqref{N2 N3 eqn defn} of our fixed point equation for $a$.
One avenue toward this study would be the oscillatory integral theory developed by Lombardi \cite{lombardi} in analytic function spaces and invoked in a spatial dynamics approach to several FPUT lattice problems \cite{IK, iooss, venney-zimmer}.
We are optimistic that techniques developed for further ripple analysis in any one of these heterogeneous lattice problems, FPUT or MiM, would readily translate to all the others due to the fundamentally similar structure of their traveling wave problems.

Last, we remark that the strictly linear bead-resonator force $r \mapsto \kappa{r}$ was very convenient for our construction and analysis of the traveling wave problem; for this reason, it is a popular force in the literature.
In the antiresonance case, it facilitates in both \cite{ksx} and \cite{faver-goodman-wright} the reduction from two equations for two traveling wave profiles to only one equation, a key part of the analysis of each of these papers.
It would therefore be interesting to add higher order terms to the bead-resonator force and see if the resulting complications are intrinsically challenging or merely technical.

\appendix

\section{Fourier Analysis}

\subsection{Fourier transforms and series}\label{fourier appendix}
We use the following conventions for Fourier transforms and Fourier series.
For $\fb \in L^1(\R,\C^m)$, we define
\[
\ft[\fb](k)
= \hat{\fb}(k)
:= \frac{1}{\sqrt{2\pi}}\int_{-\infty}^{\infty} \fb(x)e^{-ikx} \dx
\]
and
\[
\ft^{-1}[\fb](x)
= \check{\fb}(x)
:= \frac{1}{\sqrt{2\pi}}\int_{-\infty}^{\infty} \fb(k)e^{ikx} \dx.
\]
If $\fb \in L_{\per}^2(\R,\C^m)$, we define 
\[
\ft[\fb](k)
= \hat{\fb}(k)
= \frac{1}{2\pi}\int_{-\pi}^{\pi} \fb(x)e^{-ik{x}} \dx,
\]
and we have
\[
\fb(x)
= \sum_{k=-\infty}^{\infty} \hat{\fb}(k)e^{ik{x}}.
\]
We regularly use the following identity governing the interaction of Fourier transforms and series with the shift operator $S^d$:
\[
\hat{S^d\fb}(k) = e^{ikd}\hat{\fb}(k).
\]

\subsection{Fourier multipliers}

Suppose that $\tM \colon \R \to \C^{m \times m}$ is measurable and $r$, $s \in \R$ with 
\[
\sup_{k \in \R} \frac{\norm{\tM(k)}}{(1+k^2)^{(r-s)/2}}
< \infty.
\]
Then we may define the Fourier multiplier $\M \colon H^r(\R,\C^m) \to H^s(\R,\C^m)$ by 
\[
\hat{\M\fb}(k) = \tM(k)\hat{\fb}(k).
\]

The following is a slight generalization of Lemma 3 in \cite{beale1}, phrased in \cite{faver-hupkes} as Lemma A.3.
For $0 \le q \le \grave{q}$, we set
\[
S_{q,\grave{s}}
:= \set{z \in \C}{q \le |z| \le \grave{q}}.
\]

\begin{lemma}[Beale]\label{beale fm}
Let $0 < q_0 \le q_1 < q_2$ and suppose that $\tM \colon \R \to \C$ is a measurable function with the following properties.

\begin{enumerate}[label={\bf($\M$\arabic*)}]

\item
The function $\tM$ is analytic on the strips $S_{0,q_1}$ and $S_{q_1,q_2}$.

\item\label{beale simple zeros}
The function $\tM$ has finitely many zeros in $\R$, all of which are simple.  
Denote the collection of these zeros by $\Pfrak_{\M}$.

\item
There exist $C$, $z_0 > 0$ and $s \ge 0$ such that if $z \in S_{0,q_1} \cup S_{q_1,q_2}$ with $|z| \ge z_0$, then
\begin{equation}\label{beale est}
C|\re(z)|^s
\le |\tM(z)|.
\end{equation}
\end{enumerate}

Now let $\M$ be the Fourier multiplier with symbol $\tM$.
There exist $q_{\star}$, $q_{\star\star} > 0$ with $q_1 \le q_{\star} < q_{\star\star} \le q_2$ such that if $q \in [q_{\star},q_{\star\star}]$, then, for any $r \ge 0$, $\M$ is invertible from $H_q^{r+s}$ to the subspace
\[
\Dfrak_{\M,q}^r
:= \set{f \in H_q^r}{z \in \Pfrak_{\M} \Longrightarrow \hat{f}(z) = 0}
\]
and, for $f \in \Dfrak_{\M,q}^r$,
\[
\norm{\M^{-1}f}_{r+s,q}
\le \left(\sup_{k \in \R} \frac{(1+k^2)^{s/2}}{|\tM(k\pm{iq})|}\right)\norm{f}_{r,s}.
\]

In particular, if $q_1 = q_2$, then the result above is true for all $0 < q \le q_2$.
\end{lemma}
\section{Construction of Periodic Solutions}\label{periodic solutions appendix}

\subsection{Linear analysis}
All of the information that we need for our periodic bifurcation argument is contained in the following proposition.
We will need to recall the definitions and notation for the periodic Sobolev spaces introduced at the start of Section \ref{periodic solutions section}.

\begin{proposition}\label{eigenprops}
The following statements are true when $\mu$ is sufficiently small.

\begin{enumerate}[label={\bf(\roman*)}, ref={(\roman*)}]

\item\label{eig bounds part}
The eigenvalues $\lambda_{\mu}^{\pm}$ are real, even, bounded $2\pi$-periodic functions.
In particular,
\[
0 
\le \lambda_{\mu}^-(K)
\le \lambda_{\mu}^-(\pi)
< \kappa(1+\mu)
= \lambda_{\mu}^+(0)
< \lambda_{\mu}^+(K)
\le \lambda_{\mu}^+(\pi).
\]

\item\label{-zeros part}
$c^2\mu{K}^2-\lambda_{\mu}^-(K) = 0$ if and only if $K = 0$.

\item\label{omega-mu part}
There exists a unique $\omega_{\mu} > 0$ such that 
\[
c^2\mu\omega_{\mu}^2-\lambda_{\mu}^+(\omega_{\mu}) = 0.
\]
Moreover,
\begin{equation}\label{omega-mu estimates}
\Omega_{\mu}
= \frac{\sqrt{\kappa(1+\mu)}}{c\sqrt{\mu}} 
\le \omega_{\mu}
\le \frac{\sqrt{\lambda_{\mu}^+(\pi)}}{c\sqrt{\mu}}
\end{equation}
and
\begin{equation}\label{omega-mu vs Omega-mu}
\omega_{\mu}-\Omega_{\mu}
= \O(\mu^{3/2}).
\end{equation}

\item\label{upsilon-mu part}
There exists $\alpha_{\mu} \in \R$ with $\alpha_{\mu} = \O(1)$ such that if 
\begin{equation}\label{upsilon-mu defn}
\upsilon_{\mu} := \sin\left(\frac{\omega_{\mu}}{2}\right)\alpha_{\mu}
\end{equation}
and
\begin{equation}\label{nub-mu defn}
\nub_{\mu}(X) := \begin{pmatrix*}
\upsilon_{\mu}\cos(X) \\ \sin(X) 
\end{pmatrix*},
\end{equation}
then the kernel of $\Gamma_{\mu}[\omega_{\mu}]$ in $E_{\per,0}^{r+2} \times O_{\per}^{r+2}$, $r \ge 0$, is spanned by $\nub_{\mu}$.

\item\label{adj kernel part}
There exists $z_{\mu} \in \R$ with $z_{\mu} = \O(\mu^2)$ such that if 
\begin{equation}\label{nub-mu-star defn}
\nub_{\mu}^*(X) := \begin{pmatrix*}
z_{\mu}\cos(X) \\ \sin(X)
\end{pmatrix*},
\end{equation}
then, given $\fb \in E_{\per,0}^{r+2} \times O_{\per}^{r+2}$ and $\gb \in E_{\per,0}^r \times O_{\per}^r$, we have $\Gamma_{\mu}[\omega_{\mu}]\fb = \gb$ if and only if $\ip{\gb}{\nub_{\mu}^*}_0 = 0$.

\item\label{coercive part}
There is a constant $C_r > 0$, independent of $\mu$, such that if $\Gamma_{\mu}[\omega_{\mu}]\fb = \gb$ for $\fb \in E_{\per,0}^{r+2} \times O_{\per}^{r+2}$ and $\gb \in E_{\per,0}^r \times O_{\per}^r$ with $\ip{\fb}{\nub_{\mu}^*}_0 = \ip{\gb}{\nub_{\mu}^*}_0 = 0$, then $\norm{\fb}_{r+2} \le C_r\norm{\gb}_r$.
\end{enumerate}
\end{proposition}

\begin{proof}
\begin{enumerate}[label={\bf(\roman*)}]

\item
From the definition of $P_{\mu}$ in \eqref{Pmu}, we see that $P_{\mu}(\cos(K)) = \kappa^2 + \O(\mu)$, so by taking $\mu$ sufficiently small we obtain $P_{\mu}(\cos(K)) > 0$, and therefore $\lambda_{\pm}^{\mu}(K)$ are real for all $K$.
That $\lambda_{\mu}^{\pm}$ are even, $2\pi$-periodic, and continuously differentiable is then obvious.

For the bounds, we first show that the critical points of $\lambda_{\mu}^{\pm}$ are precisely the set $\pi\Z$.
We compute
\begin{equation}\label{lambda deriv formula}
(\lambda_{\mu}^{\pm})'(K)
= \mu\sin(K)\left(1 \mp \frac{2\mu\cos(K)+\kappa-(\kappa+2)\mu}{\sqrt{P_{\mu}(\cos(K))}}\right).
\end{equation}
Next, we claim that when $\mu$ is sufficiently small, we have
\begin{equation}\label{lambda deriv ineq intermediate}
0 < \frac{2\mu\cos(K)+\kappa-(\kappa+2)\mu}{\sqrt{P_{\mu}(\cos(K))}} < 1.
\end{equation}
If this is true, then the only zeros of $(\lambda_{\mu}^{\pm})'$ will be those of $\sin(\cdot)$.
It is straightforward to check that if 
\[
\mu < \min\left\{\frac{\kappa}{8},\frac{\kappa}{4(\kappa+2)}\right\}, 
\]
then we have 
\[
2\mu\cos(K)+\kappa-(\kappa+2)\mu > 0
\]
for all $K$.
After some elementary algebraic rearrangements, we see that the inequality 
\[
2\mu\cos(K)+\kappa-(\kappa+2)\mu 
< \sqrt{P_{\mu}(\cos(K))}
\]
is equivalent to $4\kappa^2\mu > 0$, which is always true.
And so \eqref{lambda deriv ineq intermediate} holds, and therefore the critical points of $(\lambda_{\mu})^{\pm}$ are $\pi\Z$.

Corollary to \eqref{lambda deriv formula} and \eqref{lambda deriv ineq intermediate} are the useful inequalities
\begin{equation}\label{lambda deriv ineqs}
|(\lambda_{\mu}^+)'(K)| \le \mu|K|
\quadword{and}
|(\lambda_{\mu}^-)'(K)| \le 2\mu|K|,
\end{equation}
which we will need in later parts of the proof.

Now we are able to bound $\lambda_{\mu}^{\pm}$.
We refer to \eqref{lambda deriv formula} to see that $(\lambda_{\mu}^{\pm})'(K) = 0$ if and only if $\sin(K) = 0$, which happens if and only if $K \in \pi\Z$.
Since $\lambda_{\mu}^{\pm}$ are $2\pi$-periodic, we only need to evaluate these functions at $K = 0$ and $K=\pi$.
We compute
\[
\lambda_{\mu}^-(0) 
= \frac{\kappa(1+\mu)}{2} - \frac{\sqrt{P_{\mu}(1)}}{2}
= \frac{\kappa(1+\mu)}{2} - \frac{\sqrt{\kappa^2(1+\mu)^2}}{2}
= 0
\]
and
\[
\lambda_{\mu}^-(\pi)
= \frac{\kappa(1+\mu)}{2} +2\mu- \frac{\sqrt{P_{\mu}(-1)}}{2}.
\]
The quantity $P_{\mu}(-1)$ will also be important for understanding $\lambda_{\mu}^+(\pi)$, so we provide the estimates
\begin{equation}\label{Pmu(-1) estimates}
\frac{\sqrt{P_{\mu}(-1)}}{2} < \frac{\kappa(1+\mu)}{2}+2\mu
\quadword{and}
\frac{\kappa(1+\mu)}{2} < 2\mu + \frac{\sqrt{P_{\mu}(-1)}}{2},
\end{equation}
which are proved with a few easy rearrangements of some equivalent inequalities.
We conclude
\[
\frac{\kappa(1+\mu)}{2} +2\mu- \frac{\sqrt{P_{\mu}(-1)}}{2} > 0,
\]
and so
\[
0 \le \lambda_{\mu}^-(K) \le \lambda_{\mu}^-(\pi) < \kappa(1+\mu).
\]
for all $K$.

Next, we evaluate
\[
\lambda_{\mu}^+(0)
= \kappa(1+\mu)
\quadword{and}
\lambda_{\mu}^+(\pi)
= \frac{\kappa(1+\mu)}{2}+2\mu + \frac{\sqrt{P_{\mu}(-1)}}{2}.
\]
By \eqref{Pmu(-1) estimates}, we have
\[
\frac{\kappa(1+\mu)}{2}+2\mu + \frac{\sqrt{P_{\mu}(-1)}}{2}
> \frac{\kappa(1+\mu)}{2} + \frac{\kappa(1+\mu)}{2}
= \kappa(1+\mu)
= \lambda_{\mu}^+(0), 
\]
and so
\[
\lambda_{\mu}^+(0) 
= \kappa(1+\mu)
\le \lambda_{\mu}^+(K)
\le \lambda_{\mu}^+(\pi)
\]
for all $K$.

\item
We use the second inequality in \eqref{lambda deriv ineqs} and the calculation $\lambda_{\mu}^-(0) = 0$ to estimate
\begin{multline*}
c^2\mu{K}^2-\lambda_{\mu}^-(K)
= \int_0^K (2c^2\mu{s}-(\lambda_{\mu}^-)'(s)) \ds
\ge \int_0^K (2c^2\mu{s}-2\mu{s}) \ds
= \int_0^K 2\mu(c^2-1)s \ds \\
\\
= \mu(c^2-1)K^2.
\end{multline*}
Since $c^2 > 1$, we have
\[
c^2K^2-\lambda_{\mu}^-(K) > 0 
\]
for all $K > 0$, and so by evenness $c^2K^2-\lambda_{\mu}^-(K) \ne 0$ for all $K \ne 0$.

\item
Since
\[
\kappa(1+\mu) \le \lambda_{\mu}^+(K) \le \lambda_{\mu}^+(\pi)
\]
for all $K$, we see that if $0 < K < \sqrt{\kappa(1+\mu)/c^2\mu} = \Omega_{\mu}$, then
\[
c^2\mu{K}^2 - \lambda_{\mu}^+(K) < 0,
\]
while if $K > \sqrt{\lambda_{\mu}^+(\pi)/c^2\mu}$, then
\[
c^2\mu{K}^2 - \lambda_{\mu}^+(K) > 0.
\]
Consequently, there exists at least one $\omega_{\mu} > 0$ with
\[
\Omega_{\mu}
\le \omega_{\mu}
\le \sqrt{\frac{\lambda_{\mu}^+(\pi)}{c^2\mu}}
\]
such that 
\[
c^2\mu\omega_{\mu}^2 - \lambda_{\mu}^+(\omega_{\mu}) = 0.
\]

For the uniqueness of $\omega_{\mu}$, we use the first inequality in \eqref{lambda deriv ineqs} to bound
\[
\frac{d}{dK}[c^2\mu{K}^2-\lambda_{\mu}^+(K)]
= 2c^2\mu{K}-(\lambda_{\mu}^+)'(K)
\ge (2c^2-1)\mu{K}
> 0
\]
when $K > 0$, so $c^2\mu{K}^2-\lambda_{\mu}^+(K)$ is strictly increasing on $(0,\infty)$ and therefore has at most one positive root.
By evenness this root is unique on all of $\R$.

Last, we prove the estimate \eqref{omega-mu vs Omega-mu}.
We expand the square root in the definition of $\lambda_{\mu}^+$ with Taylor's theorem to calculate
\begin{multline}\label{omega-mu O(mu2)}
0 = c^2\mu\omega_{\mu}^2 - \lambda_{\mu}^+(\omega_{\mu})
= c^2\mu\omega_{\mu}^2 - \frac{\kappa(1+\mu)}{2} + \frac{(2\cos(\omega_{\mu})-2)\mu}{2} \\
\\ 
- \frac{1}{2}\left(\sqrt{\kappa^2} + \frac{1}{2\sqrt{\kappa^2}}(2\kappa)(2\cos(\omega_{\mu})-2+\kappa)\mu + \O(\mu^2)\right) 
= c^2\mu\omega_{\mu}^2 - \kappa(1+\mu) + \O(\mu^2).
\end{multline}
Then, recalling that
\[
c^2\mu\Omega_{\mu}-\kappa(1+\mu) 
= 0,
\]
we subtract this equation from \eqref{omega-mu O(mu2)} to find
\[
c^2\mu(\omega_{\mu}+\Omega_{\mu})(\omega_{\mu}-\Omega_{\mu})
= \O(\mu^2).
\]
Since $\omega_{\mu}$ and $\Omega_{\mu}$ are both $\O(\mu^{-1/2})$, we conclude that
\[
\omega_{\mu}-\Omega_{\mu}
= \O(\mu^{3/2}).
\]

\item
If $\Gamma_{\mu}[\omega_{\mu}]\phib = 0$ and $\hat{\phib}(k) \ne 0$, then either
\begin{equation}\label{kernel pos}
c^2\mu(\omega_{\mu}k)^2 = \lambda_{\mu}^-(\omega_{\mu}k)
\quadword{or}
c^2\mu(\omega_{\mu}k)^2 = \lambda_{\mu}^+(\omega_{\mu}k).
\end{equation}
By part \ref{-zeros part}, the first case holds if and only if $k =0$.
But we have assumed that $\phi_1$ is mean-zero and $\phi_2$ is odd, which forces $\hat{\phib}(0) = 0$, so only the second equality in \eqref{kernel pos} is possible.

In that case, part \ref{omega-mu part} implies $\omega_{\mu}k = \pm\omega_{\mu}$, so $k=\pm1$.
We conclude that if $\Gamma_{\mu}[\omega_{\mu}]\phib = 0$ and $\hat{\phib}(k) \ne 0$, then $k=\pm1$, so
\begin{equation}\label{initial phib}
\phib(X) = e^{-iX}\hat{\phib}(-1) + e^{iX}\hat{\phib}(1),
\end{equation}
where $\hat{\phib}(1)$ is an eigenvector of $\tD_{\mu}[\omega_{\mu}]$ corresponding to the eigenvalue $\lambda_{\mu}^+(\omega_{\mu})$.
Since $\phi_1$ is even and $\phi_2$ is odd, \eqref{initial phib} rearranges into 
\begin{equation}\label{phib2}
\phib(X) = 2\begin{pmatrix*}
\hat{\phi}_1(1)\cos(X) \\
i\hat{\phi}_2(1)\sin(X)
\end{pmatrix*}.
\end{equation}

The definition of $\tD_{\mu}$ in \eqref{t-D-mu defn} and this eigenvalue-eigenvector relationship imply 
\begin{equation}\label{alpha-mu defn}
\hat{\phib}(1) = a\begin{pmatrix*}
i\alpha_{\mu}\sin\left(\frac{\omega_{\mu}}{2}\right) \\
1
\end{pmatrix*},
\qquad
\alpha_{\mu} := \frac{2\kappa}{\mu(2\cos(\omega_{\mu})-2)+\lambda_{\mu}^+(\omega_{\mu})}
\end{equation}
for some scalar $a$.
Note that the definition of $\lambda_{\mu}^+$ in \eqref{lambda defn} gives
\[
(2\cos(\omega_{\mu})-2)\mu+\lambda_{\mu}^{\pm}(\omega_{\mu})
= \frac{\kappa}{2} + \frac{\sqrt{P_{\mu}(\cos(\omega_{\mu}))}}{2} + \O(\mu),
\]
and so for $\mu$ sufficiently small the denominator of $\alpha_{\mu}$ is bounded away from zero.
That is, $\alpha_{\mu} = \O(1)$, and combining \eqref{phib2} and \eqref{upsilon-mu defn} shows that $\phib = a\nub_{\mu}$ for some scalar $a$, where $\nub_{\mu}$ is defined as in \eqref{nub-mu defn}.

\item
It suffices to show that the kernel of $\Gamma_{\mu}[\omega_{\mu}]^*$ is spanned by a vector $\nub_{\mu}^*$ of the form given in \eqref{nub-mu-star defn}.
Here $\Gamma_{\mu}[\omega_{\mu}]^* \colon E_{\per,0}^r \times O_{\per}^r \to E_{\per,0}^{r+2} \times O_{\per}^{r+2}$ is the adjoint of $\Gamma_{\mu}[\omega_{\mu}]$.
Then part \ref{adj kernel part} will follow by classical functional analysis.

Elementary theory of Fourier multipliers tells us that this adjoint satisfies
\[
\ft[\Gamma_{\mu}[\omega_{\mu}]^*\phib](k)
= \frac{1}{(1+k^2)^2}(-c^2\mu\omega_{\mu}^2k^2 + I^{\mu}\tD_{\mu}[\omega_{\mu}k]^*)\hat{\phi}(k),
\]
where $\tD_{\mu}[K]^*$ is the conjugate transpose of the matrix $\tD_{\mu}[K]$.
Now let $\phib \in E_{\per,0}^r \times O_{\per}^r$ be a nontrivial element of the kernel of $\Gamma_{\mu}[\omega_{\mu}]^*$.
Then $c^2\mu\omega_{\mu}k^2$ is an eigenvalue of $I^{\mu}\tD_{\mu}[\omega_{\mu}k]^*$ corresponding to the eigenvector $\hat{\phib}(k)$.
Since $I^{\mu}\tD_{\mu}[\omega_{\mu}k]^* = (\tD_{\mu}[\omega_{\mu}k]I^{\mu})^*$, and the eigenvalues of $\tD_{\mu}[\omega_{\mu}k]I^{\mu}$ are the real-valued functions $\lambda_{\mu}^{\pm}(\omega_{\mu}k)$, the same critical alternative \eqref{kernel pos} must hold.
Repeating the analysis from part \ref{upsilon-mu part}, we obtain \eqref{phib2} again, so we calculate $\hat{\phib}(1)$.
This time we find
\[
\hat{\phib}(1)
= a\begin{pmatrix*}
z_{\mu} \\
1
\end{pmatrix*},
\qquad
z_{\mu} := \frac{2\mu^2}{\mu(2\cos(\omega_{\mu})-2)+\lambda_{\mu}^+(\omega_{\mu})}\sin\left(\frac{\omega_{\mu}}{2}\right),
\]
so $z_{\mu} = \O(\mu^2)$.

\item 
If $\Gamma_{\mu}[\omega_{\mu}]\fb = \gb$ for some $\fb \in E_{\per,0}^{r+2}\times O_{\per}^{r+2}$ and $\gb \in E_{\per,0}^r \times O_{\per}^r$, then, for $k \ne 0$, we have
\begin{equation}\label{pre-neumann}
\left(\ind - \frac{1}{c^2\mu\omega_{\mu}^2k^2}\tD_{\mu}[\omega_{\mu}k]I^{\mu}\right)\hat{\fb}(k) 
= -\frac{1}{c^2\mu\omega_{\mu}^2k^2}\hat{\gb}(k).
\end{equation}
Here $\ind$ is the $2\times2$ identity matrix.
For $|k| \ge 2$, we will show that the $\infty$-norm of the matrix
\[
\frac{1}{c^2\omega_{\mu}^2k^2}\tD_{\mu}[\omega_{\mu}k]I^{\mu}
\]
is smaller than 1 uniformly in $\mu$ and $k$.
This, along with estimating $\mu\omega_{\mu}^2 = \O(1)$ on the right side of \eqref{pre-neumann}, enables us to invert the matrix on the left side and obtain an estimate of the form $|\hat{\fb}(k)| \le C|\hat{\gb}(k)|/k^2$ for $|k| \ge 2$.
When $|k| = 1$, we will prove directly that $|\hat{\fb}(\pm1)| \le C|\hat{\fb}(\pm1)|$.
Combining these estimates, Plancherel's theorem will imply the desired inequality $\norm{\fb}_{r+2} \le C_r\norm{\gb}_r$.

Referring to the definition of $\tD_{\mu}$ in \eqref{t-D-mu defn}, we see that by taking $\mu$ sufficiently small, we have
\[
\norm{\tD_{\mu}[\omega_{\mu}k]I^{\mu}}_{\infty}
\le 2\kappa.
\]
By the estimates on $\omega_{\mu}$ from \eqref{omega-mu estimates}, this implies
\[
\frac{1}{c^2\mu\omega_{\mu}^2k^2}\norm{\tD_{\mu}[\omega_{\mu}k]I^{\mu}}_{\infty}
\le \frac{2}{k^2} 
< \frac{1}{2}
\]
when $|k| \ge 2$.
A Neumann series argument allows us to conclude
\[
\bignorm{\left(\ind - \frac{1}{c^2\mu\omega_{\mu}^2k^2}\tD_{\mu}[\omega_{\mu}k]I^{\mu}\right)^{-1}}_{\infty} 
\le 2.
\]

Now we consider the case $|k| = 1$.
In this case, the equality $\hat{\Gamma_{\mu}[\omega_{\mu}]\fb}(1)\cdot\e_1 = \hat{\gb}(1)\cdot\e_1$ is equivalent to
\begin{equation}\label{matrix inverse1}
-c^2\mu\omega_{\mu}^2\hat{f}_1(1) -\mu(2\cos(\omega_{\mu})-2)\hat{f}_1(1)+2i\kappa\sin\left(\frac{\omega_{\mu}}{2}\right)\hat{f}_2(1) 
= \hat{g}_1(1).
\end{equation}
It is here that we need the additional hypothesis from part \ref{coercive part} that 
\[
\ip{\fb}{\nub_{\mu}^*}_0 
= \ip{\gb}{\nub_{\mu}^*}_0
= 0.
\]
Since $f_1$ and $g_1$ are even and $f_2$ and $g_2$ are odd, these equalities and the definition of $\nub_{\mu}^*$ in \eqref{nub-mu-star defn} imply
\[
\hat{f}_2(1) = -iz_{\mu}\hat{f}_1(1)
\quadword{and}
\hat{g}_2(1) = -iz_{\mu}\hat{g}_1(1).
\]
We can therefore remove $\hat{f}_2(1)$ from the left side of \eqref{matrix inverse1} and rewrite that equation as
\begin{equation}\label{matrix inverse2}
\left(-c^2\mu\omega_{\mu}^2 - \mu(2\cos(\omega_{\mu})-2) + 2z_{\mu}\kappa\sin\left(\frac{\omega_{\mu}}{2}\right)\right)\hat{f}_1(1) 
= \hat{g}_1(1).
\end{equation}
Since $\mu\omega_{\mu}^2 = \O(1)$ and $z_{\mu} = \O(\mu^2)$, the prefactor on $\hat{f}_1(1)$ is bounded away from zero uniformly in $\mu$ when $\mu$ is small.
We conclude $|\hat{f}_1(1)| \le C|\hat{g}_1(1)|$.

If $z_{\mu} = 0$, then $|\hat{f}_2(1)| = 0 = |\hat{g}_2(1)|$, and there is nothing to prove for the estimate on the second component.
Otherwise, we can write
\[
\frac{1}{|z_{\mu}|} |\hat{f}_2(1)|
= |\hat{f}_1(1)|
\le C|\hat{g}_1(1)|
= \frac{C}{|z_{\mu}|}|\hat{g}_2(1)|,
\]
and thus $|\hat{f}_2(1)| \le C|\hat{g}_2(1)|$.

The estimates for $k=-1$ now follow by the even-odd symmetry of $\fb$ and $\gb$.
\qedhere
\end{enumerate}
\end{proof}

\subsection{The bifurcation argument}
Proposition \ref{eigenprops} is analogous to Lemma B.1 in \cite{hoffman-wright}, which contained all the information necessary for their construction of the periodic traveling waves.
Since the argument for the MiM periodic solutions would proceed exactly as the one for the small mass FPUT periodics, we do not feel the need to give it here and sketch it only fleetingly; a longer summary can be found in Appendix C of \cite{faver-hupkes}, and the full argument takes place in Appendix B of \cite{hoffman-wright}.
This argument bears some strong similarities with the periodic proofs for the long wave problem in Appendix C  of \cite{faver-wright} and relies on a number of Lipschitz estimates for scaled Fourier multipliers, which are detailed in \cite{faver-wright} and also Appendix D of \cite{faver-dissertation}.

Here, then, in brief is the argument.
One wishes to solve $\Phib_{\mu}(\phib,\mu) = 0$ from \eqref{periodic problem} and makes the revised ansatz
\[
\phib
= a\nub_c + a\psib
\quadword{and}
\omega
= \omega_{\mu} + \xi,
\]
where $\psib \in E_{\per,0}^2 \times O_{\per}^2$ and $\xi \in \R$.
This is similar to the ansatz made in the actual proof of the bifurcation from simple eigenvalues theorem by Crandall and Rabinowitz \cite{crandall-rabinowitz}.
Then one uses the various properties delineated in Proposition \ref{eigenprops} to rearrange the problem $\Phib_{\mu}(a\nub_c+a\psib,\omega_{\mu}+\xi)$ into a fixed point problem for $\psib$ and $\xi$, where $a$ and $\mu$ are parameters.
Specifically, one finds an equation for $\xi$ by projecting onto the one-dimensional span of $\nub_c$ and then an equation for $\psib$ by projecting onto the orthogonal complement of this span; parts \ref{adj kernel part} and \ref{coercive part} of Proposition \ref{eigenprops} are critical here.
A quantitative fixed point argument, using many of the Fourier multiplier estimates discussed above, finishes the job.
\section{Proofs for the Construction of the Nanopteron Fixed Point Argument}

\subsection{Proof of Proposition \ref{L prop}}\label{L prop proof appendix}
This proof is similar to the proofs of Lemmas A.8 and A.9 in \cite{faver-wright}.
However, the notation and estimates are sufficiently different that we feel it worthwhile to spell out the details here.

Fix $|c| > 1$, $\kappa > 0$, and $\mu \in (0,1)$, and let
\[
\tM(z)
:= c^2z^2-\kappa(1+\mu),
\]
so that 
\[
\hat{\L_{\mu}f}(k) = \tM(\sqrt{\mu}k)\hat{f}(k).
\]
Since $\tM$ is just a quadratic polynomial with the distinct roots $\pm\sqrt{\kappa(1+\mu)}/|c|$, it is straightforward to see that all the hypotheses of Lemma \ref{beale fm} are satisfied.
In particular, they hold for any $q > 0$, since the only zeros of $\tM$ are real.

Lemma \ref{beale fm} also gives the estimate
\[
\norm{\L_{\mu}^{-1}}_{\b(\Dfrak_{\mu,q}^r,O_q^{r+j})}
\le \sup_{k \in \R} \frac{(1+k^2)^{j/2}}{|\tM(\sqrt{\mu}(k+iq))|}.
\]
We estimate this supremum further:
\begin{multline*}
\sup_{k \in \R} \frac{(1+k^2)^{j/2}}{|\tM(\sqrt{\mu}(k+iq))|}
= \sup_{k \in \R} \frac{(1+(k/\sqrt{\mu})^2)^{j/2}}{|\tM(k+i\sqrt{\mu}q)|}
= \sup_{k \in \R} \left(\frac{(1+(k/\sqrt{\mu})^2)^{j/2}}{(1+k^2)^{j/2}}\right)\left(\frac{(1+k^2)^{j/2}}{\tM(k+i\sqrt{\mu}q)}\right) \\
\\
\le \left(\sup_{k \in \R}\frac{(1+(k/\sqrt{\mu})^2)^{j/2}}{(1+k^2)^{j/2}} \right)\left(\sup_{k \in \R}\frac{(1+k^2)^{j/2}}{\tM(k+i\sqrt{\mu}q)} \right).
\end{multline*}
Some elementary algebraic rearrangements imply
\[
\sup_{k \in \R}\frac{(1+(k/\sqrt{\mu})^2)^{j/2}}{(1+k^2)^{j/2}}
\le \frac{1}{\mu^{j/2}}.
\]
We claim more generally that if $0 < |\tau| < \tau_0$ for some $\tau_0 > 0$, then there is a constant $C > 0$, independent of $\tau$ (but dependent on $c$, $\kappa$, and $\tau_0$) such that 
\begin{equation}\label{tau claim}
\sup_{k \in \R}\frac{(1+k^2)^{j/2}}{\tM(k+i\tau)}
\le \frac{C}{|\tau|}
\end{equation}
for $j=0,1,2$. 
From this our desired estimate \eqref{L-mu coercive} follows.

So, we just need to prove this claim \eqref{tau claim}.
First, if $z \in \C$ with $|z| \ge 2\sqrt{\kappa}/|c|$, then the reverse triangle inequality implies
\[
|\tM(z)| 
\ge \frac{c^2|z|^2}{2},
\]
from which we see that if $k \in \R$ with $|k| \ge 2\sqrt{\kappa}/|c|$ and $0 < |\tau| < \tau_0$, then 
\[
\frac{|\tM(k+i\tau)|}{1+k^2}
\ge \frac{c^2|k+i\tau|^2}{2(1+k^2)}
\ge \frac{c^2k^2}{2(1+k^2)}
\ge \frac{c^2}{4}
= \left(\frac{c^2}{4\tau_0}\right)\tau_0
\ge \left(\frac{c^2}{4\tau_0}\right)|\tau|.
\]
Then for $j=0,1,2$,
\[
\sup_{\substack{k \in \R \\ |k| \ge 2\sqrt{\kappa}/|c|}} \frac{(1+k^2)^j/2}{|\tM(k+i\tau)|}
\le \sup_{\substack{k \in \R \\ |k| \ge 2\sqrt{\kappa}/|c|}} \frac{1+k^2}{|\tM(k+i\tau)|}
\le \left(\frac{4\tau_0}{c^2}\right)\frac{1}{|\tau|}.
\]
When $|k| \le 2\sqrt{\kappa}/|c|$, we use a partial fractions decomposition to rewrite
\[
\frac{1}{\tM(z)}
= -\frac{1}{2\sqrt{\kappa(1+\mu)}\big(cz+\sqrt{\kappa(1+\mu)}\big)}
+ \frac{1}{2\sqrt{\kappa(1+\mu)}\big(cz-\sqrt{\kappa(1+\mu)}\big)},
\]
and then it is apparent that 
\[
\frac{1}{|c(k+i\tau)\pm\sqrt{\kappa(1+\mu)}|}
= \frac{1}{\sqrt{\big(ck\pm\sqrt{\kappa(1+\mu)}\big)^2+c^2\tau^2}}
\le \frac{1}{|c\tau|}.
\]
Thus
\[
\sup_{\substack{k \in \R \\ |k| \le 2\sqrt{\kappa}/|c|}} \frac{1+k^2}{|\tM(k+i\tau)|}
\le \left(\frac{1+4\kappa/c^2}{|c|}\right)\frac{1}{|\tau|}.
\]
This completes the proof of \eqref{tau claim} and also of Proposition \ref{L prop}.

\subsection{Proof of the identity \eqref{chi-id}}\label{chi-id proof appendix}
With $\chi_{\mu}$ defined in \eqref{chi-mu defn} and $\iota_{\mu}$ in \eqref{iota-mu defn}, we calculate
\begin{multline*}
\iota_{\mu}[\chi_{\mu}]
= -2\mu^2\upsilon_{\mu}\int_{-\infty}^{\infty} \big(\delta\big[\sigma_c\cos(\Omega_{\mu}\cdot)\big]\big)(x)\sin(\Omega_{\mu}x) \dx \\
\\
= 2\mu^2\upsilon_{\mu}\int_{-\infty}^{\infty} \sigma_c(x)\cos(\Omega_{\mu}x)\big(\delta\big[\sin(\Omega_{\mu}\cdot)\big]\big)(x) \dx
= 4\mu^2\upsilon_{\mu}\sin\left(\frac{\Omega_{\mu}}{2}\right)\int_{-\infty}^{\infty} \sigma_c(x)\cos^2(\Omega_{\mu}x) \dx.
\end{multline*}
To obtain the second equality, we used the skew-adjointness of $\delta$ in the $L^2$-inner product, and for the third we directly calculated $\delta[\sin(\Omega_{\mu}\cdot)]$ using the addition formulas for sine.

Next, we have
\[
\sigma_c(x)\cos^2(\omega_{\mu}x) 
= \sigma_c(x)\left(\frac{e^{i\omega_{\mu}x}+e^{-i\omega_{\mu}x}}{2}\right)^2
= \frac{\sigma_c(x)}{4} + \frac{\sigma_c(x)}{4}(e^{2i\omega_{\mu}x}+e^{-2i\omega_{\mu}x}).
\]
A variation on the Riemann-Lebesgue lemma (such as that given in Lemma A.4 of \cite{faver-wright}, which also proves Proposition \ref{iota-mu prop} in this paper) shows
\[
\int_{-\infty}^{\infty} \frac{\sigma_c(x)}{4}(e^{2i\omega_{\mu}x}+e^{-2i\omega_{\mu}x}) \dx 
= \O(\mu^{\infty}).
\]
Since $\sigma_c(x) > 0$ for all $x$, we conclude
Thus
\[
\iota_{\mu}[\chi_{\mu}] 
= \mu^2\upsilon_{\mu}\sin\left(\frac{\omega_{\mu}}{2}\right)\norm{\sigma_c}_{L^1} + \O(\mu^{\infty}).
\]

Now, recall from \eqref{upsilon pseudo-defn} that 
\[
\upsilon_{\mu}
= \sin\left(\frac{\omega_{\mu}}{2}\right)\alpha_{\mu},
\]
where $\alpha_{\mu} = \O(1)$.
Also, part \ref{omega vs Omega} of Lemma \ref{omega-mu lemma} tells us 
\[
\omega_{\mu}-\Omega_{\mu}
= \O(\mu^{3/2}).
\]
It follows that
\[
\iota_{\mu}[\chi_{\mu}] = \mu^2\sin^2\left(\frac{\Omega_{\mu}}{2}\right)\alpha_{\mu}\norm{\sigma_c}_{L^1} + \O(\mu^{7/2}).
\]
\section{The Proof of Theorem \ref{main nanopteron theorem}}

For simplicity in this appendix, we abbreviate $q_{\star} := q_{\star}(c)$ and $\qbar_{\star}(c) := \qbar_{\star}$.
These decay rate constants were selected in Section \ref{nano soln}.

\subsection{Proof of Lemma \ref{main workhorse lemma}}\label{proof of main workhorse lemma appendix}

We rely on the following lemma, whose proof is in Appendix \ref{proof of workhorse lemma for main workhorse lemma appendix}.

\begin{lemma}\label{workhorse lemma for main workhorse lemma}
For all $r \ge 3$, there is $C_{\star}(r) > 0$ such that the following estimates hold for any $\mu \in \Mu_{c,\kappa} \cap (0,\mu_{\chi})$.

\begin{enumerate}[label={\bf(\roman*)}, ref={(\roman*)}]

\item\label{workhorse2 map}
If $\etab \in E_{q_{\star}}^3 \times O_{q_{\star}}^3$ and $|a| \le a_{\per}$, then
\begin{equation}\label{workhorse2 map est}
\norm{\nanob_{\mu}(\etab,a)}_{\X_1}
\le C_{\star}(3)\left(\mu^{1/2} 
+ \mu^{1/2}|a| 
+ \mu^{1/2}\norm{\etab}_{3,q_{\star}} 
+ |a|^2 
+ \norm{\etab}_{3,q_{\star}}^2\right).
\end{equation}

\item\label{workhorse2 lip}
If $\etab$, $\grave{\etab} \in E_{q_{\star}}^3 \times O_{q_{\star}}^3$ and $|a|$, $|\grave{a}| \le a_{\per}$, then
\begin{multline}\label{workhorse2 lip est}
\norm{\nanob_{\mu}(\etab,a)-\nanob_{\mu}(\grave{\etab},\grave{a})}_{\X_0} \\
\\
\le C_{\star}(3)\left(\mu^{1/2} 
+ \norm{\etab}_{3,q_{\star}} 
+ \norm{\grave{\etab}}_{3,q_{\star}} 
+ |a| 
+ |\grave{a}|\right)
\left(\norm{\etab-\grave{\etab}}_{3,\qbar_{\star}} 
+ |a-\grave{a}|\right).
\end{multline}

\item\label{workhorse2 boot}
If $\etab \in E_{q_{\star}}^r \times O_{q_{\star}}^r$ and $|a| \le a_{\per}$, then
\begin{multline}\label{workhorse2 boot est1}
\norm{\nano_1^{\mu}(\etab,a)}_{r+1,q_{\star}} 
+ \norm{\nano_2^{\mu}(\etab,a)}_{r+1,q_{\star}} \\
\\
\le C_{\star}(r)\left(\mu^{1/2} 
+ \mu^{3/2-r/2}|a| 
+ \norm{\etab}_{r,q_{\star}} 
+ \mu^{1-r/2}|a|\norm{\etab}_{r,q_{\star}} 
+ \mu^{1-r/2}|a|^2 
+ \norm{\etab}_{r,q_{\star}}^2\right)
\end{multline}
and
\begin{equation}\label{workhorse2 boot est2}
|\nano_3^{\mu}(\etab,a)|
\le C_{\star}(r)\left(\mu^{r/2} + \mu^{r/2-1}\norm{\etab}_{r,q_{\star}} + \mu^{r/2-1}\norm{\etab}_{r,q_{\star}}^2 + |a|\norm{\etab}_{r,q_{\star}} + |a|^2\right).
\end{equation}
\end{enumerate}
\end{lemma}

Now we prove Lemma \ref{main workhorse lemma}.
Abbreviate $C_{\star} = C_{\star}(3)$ from the lemma above and set
\[
\mu_{\star}
= \min\left\{\frac{1}{2},
\mu_{\varsigma},
\mu_{\chi},
\frac{a_{\per}^2}{4C_{\star}^2},
\frac{1}{64C_{\star}^2(1+2C_{\star})^2},
\frac{1}{4C_{\star}^2(1+8C_{\star})^2}
\right\}.
\]
Let $\tau_{\star} = 2C_{\star}$.
Observe that if $(\etab,a) \in \U_{\mu,\tau_{\star}}^1$, where this set was defined in \eqref{Uhood}, and $0 < \mu < \mu_{\star}$, then
\[
|a| 
\le \tau_{\star}\mu^{3/2-1} 
= 2C_{\star}\mu^{1/2}
\le a_{\per},
\]
so $\nanob_{\mu}(\etab,a)$ is indeed defined.
In the following, fix $0 < \mu < \mu_{\star}$.

To prove part \ref{main workhorse map}, fix $(\etab,a) \in \U_{\mu,\tau_{\star}}^1$ and use \eqref{workhorse2 map est} to estinate
\begin{multline*}
\norm{\nano_1^{\mu}(\etab,a)}_{3,q_{\star}} 
+ \norm{\nano_2^{\mu}(\etab,a)}_{3,q_{\star}}
+ |\nano_3^{\mu}(\etab,a)| \\
\\
\le C_{\star}\left(\mu^{1/2} + \mu^{1/2}|a| + \mu^{1/2}\norm{\etab}_{3,q_{\star}} + |a|^2 + \norm{\etab}_{3,q_{\star}}^2\right) 
\le C_{\star}\mu^{1/2}(1+4C_{\star}\mu^{1/2} + 8C_{\star}^2\mu^{1/2}) \\
\\
= C_{\star}\mu^{1/2}(1 + \mu^{1/2}(4C_{\star})(1+2C_{\star})) 
< \frac{3}{2}C_{\star}\mu^{1/2} 
< \tau_{\star}\mu^{1/2}.
\end{multline*}
Consequently, $\nanob_{\mu}(\etab,a) \in \U_{\mu,\tau_{\star}}^1$ as well.

To prove part \ref{main workhorse lip}, take $(\etab,a)$, $(\grave{\etab},\grave{a}) \in \U_{\mu,\tau_{\star}}^1$ and estimate
\[
C_{\star}(\mu^{1/2} + \norm{\etab}_{3,q_{\star}} + \norm{\grave{\etab}}_{q_{\star}} + |a| + |\grave{a}|)
\le C_{\star}(\mu^{1/2} + 4\tau_{\star}\mu^{1/2}) 
= C_{\star}\mu^{1/2}(1+8C_{\star}) 
< \frac{1}{2}.
\]
Then part \ref{main workhorse lip} follows from \eqref{workhorse2 lip est}.

Last, for part \ref{main workhorse boot}, assume $(\etab,a) \in \U_{\mu,\tau_{\star}}^1 \cap \U_{\mu,\tau}^r$.
Since $(\etab,a) \in \U_{\mu,\tau_{\star}}^1$, $\nanob_{\mu}(\etab,a)$ is defined.
Now we use \eqref{workhorse2 boot est1} to bound $\norm{\nano_1^{\mu}(\etab,a)}_{r+1,q_{\star}}+ \norm{\nano_2^{\mu}(\etab,a)}_{r+1,q_{\star}}$ above by 
\begin{multline*}
C_{\star}(r)\left(\mu^{1/2} 
+ \mu^{3/2-r/2}|a|
+\norm{\etab}_{r,q} 
+ \mu^{1-r/2}|a|\norm{\etab}_{r,q} 
+ \mu^{1-r/2}|a|^2
+\norm{\etab}_{r,q}^2\right) \\
\\
\le C_{\star}(r)\left(\mu^{1/2}
+3\tau\mu^{1/2}
+2\tau^2\mu^{1/2}
+\tau^2\mu^{r/2-1}
\right) 
\le C_{\star}(r)(3\tau+2\tau^2)\mu^{1/2}.
\end{multline*}
Note that for $r \ge 3$, we have $r/2-1 \ge 1/2$, so we are indeed able to factor out a power of $\mu^{1/2}$ in the last inequality above.

Next, we want to show 
\[
|\nano_3^{\mu}(\etab,a)|
\le \ttau\mu^{(r+1)/2-1}
= \ttau\mu^{r/2-1/2}
\]
for some $\ttau > 0$.
We use \eqref{workhorse2 boot est2} to estimate 
\begin{align*}
|\nano_3^{\mu}(\etab,a)|
&\le C_{\star}(r)\left(\mu^{r/2}
+ \mu^{r/2-1}\norm{\etab}_{r,q}
+ \mu^{r/2-1}\norm{\etab}_{r,q}^2
+ |a|\norm{\etab}_{r,q} 
+ |a|^2\right) \\
\\
&\le C_{\star}(r)\left(\mu^{r/2}
+ 2\tau\mu^{r/2-1/2}
+ \tau^2\mu^{r/2}
+\tau^2\mu^{r-2}\right) \\
\\
&\le C_{\star}(r)\mu^{r/2-1/2}\left(
\mu^{1/2}
+2\tau
+\tau^2\mu^{1/2}
+\tau^2\mu^{r/2-3/2}\right).
\end{align*}

We set
\[
\ttau 
:= \max\left\{C_{\star}(r)(3\tau+2\tau^2), C_{\star}(\mu_{\star}^{1/2}
+2\tau
+\tau^2\mu_{\star}^{1/2}
+\tau^2\mu_{\star}^{r/2-3/2}\right\}.
\]
This concludes the proof of Lemma \ref{main workhorse lemma}.

\subsection{Proof of Lemma \ref{workhorse lemma for main workhorse lemma}}\label{proof of workhorse lemma for main workhorse lemma appendix}

\subsubsection{Proof of part \ref{workhorse2 map}}
Set
\[
\rhs_{\map}^{\mu}(\etab,a)
:= \mu^{1/2}
+ \mu^{1/2}|a|
+ \mu^{1/2}\norm{\etab}_{3,q_{\star}}
+|a|^2
+ \norm{\etab}_{3,q_{\star}}^2.
\]
We will show that
\begin{equation}\label{map goal}
\norm{\nano_1^{\mu}(\etab,a)}_{3,q_{\star}} + \norm{\nano_2^{\mu}(\etab,a)}_{3,q_{\star}} + |\nano_3^{\mu}(\etab,a)|
\le C\rhs_{\map}^{\mu}(\etab,a),
\end{equation}
which implies \eqref{workhorse2 map est}.

We first claim that 
\begin{equation}\label{h map}
\norm{h_k^{\mu}(\etab,a)}_{r,q_{\star}}
\le C(r)\Hfrak_{\map,r}^{\mu}(\etab,a), \ k=1,\ldots,4,
\end{equation}
and
\begin{equation}\label{ell map}
\norm{\ell_k^{\mu}(\etab,a)}_{r,q_{\star}}
\le C(r)\Lfrak_{\map,r}^{\mu}(\etab,a), \ k = 1,\ldots,6,
\end{equation}
where
\begin{equation}\label{H map}
\Hfrak_{\map,r}^{\mu}(\etab,a)
:= \mu^{1-r/2}|a| 
+ \mu\norm{\etab}_{r,q_{\star}} 
+ \mu^{1-r/2}|a|\norm{\etab}_{r,q_{\star}} 
+ \norm{\etab}_{r,q_{\star}}^2
\end{equation}
and
\begin{equation}\label{L map}
\Lfrak_{\map,r}^{\mu}(\etab,a)
:= \mu^2
+\mu^{3-r/2}|a| 
+ \mu\norm{\etab}_{r,q_{\star}} 
+ \mu\norm{\etab}_{r,q_{\star}}^2 
+\mu^{2-r/2}|a|\norm{\etab}_{r,q_{\star}}
+ \mu^{2-r/2}|a|^2.
\end{equation}
We prove these estimates in Appendix \ref{h ell ests appendix}.

Then from the definition of $\nano_2^{\mu}$ in \eqref{N2 N3 eqn defn} and the estimates on $\B_{\mu}$ in Proposition \ref{A B prop} with $j=0$, we have
\begin{multline}\label{N2 map}
\norm{\nano_2^{\mu}(\etab,a)}_{3,q_{\star}}
\le \sum_{k=1}^6\norm{\B_{\mu}\ell_k^{\mu}(\etab,a)}_{3,q_{\star}}
\le C\mu^{-1/2}\sum_{k=1}^6 \norm{\ell_k^{\mu}(\etab,a)}_{3,q_{\star}} \\
\\
\le C\mu^{-1/2}\Lfrak_{\map,3}^{\mu}(\etab,a)
\le C\rhs_{\map}^{\mu}(\etab,a).
\end{multline}
Next, we use the definition of $\nano_1^{\mu}$ in \eqref{h5} and the invertibility of $\H_c$ from Proposition \ref{H prop} to estimate
\begin{multline*}
\norm{\nano_2^{\mu}(\etab,a)}_{3,q_{\star}}
\le \sum_{k=1}^5 \norm{\H_c^{-1}h_k^{\mu}(\etab,a)}_{3,q_{\star}}
\le C\sum_{k=1}^4\norm{h_k^{\mu}(\etab,a)}_{1,q_{\star}} + C\norm{\nano_2^{\mu}(\etab,a)}_{3,q_{\star}}.
\end{multline*}
We have
\[
\norm{h_k^{\mu}(\etab,a)}_{1,q_{\star}} 
\le C\Hfrak_{\map,1}^{\mu}(\etab,a)
\le C\rhs_{\map}^{\mu}(\etab,a),
\]
and so, with the previous estimates on $\nano_2^{\mu}$, we find
\begin{equation}\label{N1 map}
\norm{\nano_1^{\mu}(\etab,a)}_{3,q_{\star}} 
\le C\rhs_{\map}^{\mu}(\etab,a).
\end{equation}

Last, we use the definition of $\nano_3^{\mu}$ from \eqref{N2 N3 eqn defn} and the estimates on $\A_{\mu}$ from Proposition \ref{A B prop} to produce the general estimate
\begin{equation}\label{N3 map r}
|\nano_3^{\mu}(\etab,a)|
\le \sum_{k=1}^6 |\A_{\mu}[\ell_k^{\mu}(\etab,a)]|
\le C(r)\mu^{r/2-2}\sum_{k=1}^6\norm{\ell_k^{\mu}(\etab,a)}_{r,q_{\star}}
\le C(r)\mu^{r/2-2}\Lfrak_{\map,r}^{\mu}(\etab,a).
\end{equation}
Taking $r=3$, we have
\[
|\nano_3^{\mu}(\etab,a)|
\le C\mu^{-1/2}\Lfrak_{\map,3}^{\mu}(\etab,a),
\]
and so, by \eqref{N2 map},
\begin{equation}\label{N3 map}
|\nano_3^{\mu}(\etab,a)|
\le C\rhs_{\map}^{\mu}(\etab,a).
\end{equation}
Combining \eqref{N2 map}, \eqref{N1 map}, and \eqref{N3 map} gives \eqref{map goal}.

\subsubsection{Proof of part \ref{workhorse2 lip}}
Set
\[
\rhs_{\lip}^{\mu}(\etab,a)
:= \mu^{1/2} 
+ \norm{\etab}_{3,\qbar_{\star}} 
+ \norm{\grave{\etab}}_{3,\qbar_{\star}}
+ |a|
+ |\grave{a}|
\]
We will show that 
\begin{multline}\label{lip goal}
\norm{\nano_1^{\mu}(\etab,a)-\nano_1^{\mu}(\grave{\etab},\grave{a})}_{3,\qbar_{\star}}
+ \norm{\nano_2^{\mu}(\etab,a)-\nano_2^{\mu}(\grave{\etab},\grave{a})}_{3,\qbar_{\star}}
+ |\nano_3^{\mu}(\etab,a)-\nano_3^{\mu}(\grave{\etab},\grave{a})| \\
\\
\le C\rhs_{\lip}^{\mu}(\etab,a)\left(\norm{\etab-\grave{\etab}}_{3,q_{\star}} + |a-\grave{a}|\right),
\end{multline}
which will imply \eqref{workhorse2 lip est}.

To do so, we will exploit a mismatch of regularity and work with $r=1$ for the $h_k^{\mu}$ terms but $r=3$ for the $\ell_k^{\mu}$ terms.
Namely, we claim that 
\[
\norm{h_k^{\mu}(\etab,a)-h_k^{\mu}(\grave{\etab},\grave{a})}_{1,\qbar_{\star}}
\le C\Hfrak_{\lip}^{\mu}(\etab,\grave{\etab},a,\grave{a})\left(\norm{\etab-\grave{\etab}}_{1,\qbar_{\star}} + |a-\grave{a}|\right), \ k = 1,\ldots,4,
\]
and
\[
\norm{\ell_k^{\mu}(\etab,a)-\ell_k^{\mu}(\grave{\etab},\grave{a})}_{3,\qbar_{\star}}
\le C\Lfrak_{\lip}^{\mu}(\etab,\grave{\etab},a,\grave{a})\left(\norm{\etab-\grave{\etab}}_{3,\qbar_{\star}} + |a-\grave{a}|\right), \ k = 1,\ldots,6,
\]
where
\[
\Hfrak_{\lip}^{\mu}(\etab,\grave{\etab},a,\grave{a})
:= \mu^{1/2} + \norm{\etab}_{1,q_{\star}} + \norm{\grave{\etab}}_{1,q_{\star}} + |a| + |\grave{a}|,
\]
while
\[
\Lfrak_{\lip}^{\mu}(\etab,\grave{\etab},a,\grave{a})
:= \mu + \mu^{1/2}\left(
\norm{\etab}_{3,q_{\star}} 
+ \norm{\grave{\etab}}_{3,q_{\star}}
+|a| 
+ |\grave{a}|\right).
\]
We prove these estimates in Appendix \ref{h ell ests appendix}, too.

Assuming this to be true, we use the definition of $\nano_2^{\mu}$ in \eqref{N2 N3 eqn defn} and the estimates on $\B_{\mu}$ from Proposition \ref{A B prop} to estimate.
\begin{multline}\label{N2 lip}
\norm{\nano_2^{\mu}(\etab,a)-\nano_2^{\mu}(\grave{\etab},\grave{a})}_{3,\qbar_{\star}}
\le \sum_{k=1}^6 \norm{\B_{\mu}\left(\ell_k^{\mu}(\etab,a)-\ell_k^{\mu}(\grave{\etab},\grave{a})\right)}_{3,\qbar_{\star}} \\
\\
\le C\mu^{-1/2} \sum_{k=1}^6\norm{\ell_k^{\mu}(\etab,a)-\ell_k^{\mu}(\grave{\etab},\grave{a})}_{3,\qbar_{\star}}
\le C\mu^{-1/2}\Lfrak_{\lip}^{\mu}(\etab,\grave{\etab},a,\grave{a})
\le C\rhs_{\lip}^{\mu}(\etab,\grave{\etab},a,\grave{a}).
\end{multline}
Next, we use the definition of $\nano_1^{\mu}$ from \eqref{h5} and the smoothing properties of $\H_c$ from Proposition \ref{H prop} to bound
\begin{multline*}
\norm{\nano_1^{\mu}(\etab,a)-\nano_1^{\mu}(\grave{\etab},\grave{a})}_{3,\qbar_{\star}}
\le \sum_{k=1}^5 \norm{\H_c^{-1}\left(h_k^{\mu}(\etab,a)-h_k^{\mu}(\grave{\etab},\grave{a})\right)}_{3,\qbar_{\star}} \\
\\
\le C\sum_{k=1}^4 \norm{h_k^{\mu}(\etab,a)-h_k^{\mu}(\grave{\etab},\grave{a})}_{1,\qbar_{\star}} + C\norm{\nano_2^{\mu}(\etab,a)-\nano_2^{\mu}(\grave{\etab},\grave{a})}_{3,q_{\star}} \\
\\
\le C\Hfrak_{\lip}^{\mu}(\etab,\grave{\etab},a,\grave{a}) + C\rhs_{\lip}^{\mu}(\etab,\grave{\etab},a,\grave{a}).
\end{multline*}
It is clear that $\Hfrak_{\lip}^{\mu}(\etab,\grave{\etab},a,\grave{a}) \le \rhs_{\lip}^{\mu}(\etab,\grave{\etab},a,\grave{a})$, and so we have
\begin{equation}\label{N1 lip}
\norm{\nano_1^{\mu}(\etab,a)-\nano_1^{\mu}(\grave{\etab},\grave{a})}_{3,\qbar_{\star}}
\le C\rhs_{\lip}^{\mu}(\etab,\grave{\etab},a,\grave{a}).
\end{equation}

Last, we use the definition of $\nano_3^{\mu}$ from \eqref{N2 N3 eqn defn} and the estimates on $\A_{\mu}$ from Proposition \ref{A B prop} to find
\begin{multline*}
|\nano_3^{\mu}(\etab,a)-\nano_3^{\mu}(\grave{\etab},\grave{a})|
\le \sum_{k=1}^6 |\A_{\mu}\left[\ell_k^{\mu}(\etab,a)-\ell_k^{\mu}(\grave{\etab},\grave{a})\right]|
\le C\mu^{-1/2}\sum_{k=1}^6 \norm{\ell_k^{\mu}(\etab,a)-\ell_k^{\mu}(\grave{\etab},\grave{a})}_{3,\qbar_{\star}},
\end{multline*}
and, recognizing this last term from \eqref{N2 lip}, we conclude that 
\begin{equation}\label{N3 lip}
|\nano_3^{\mu}(\etab,a)-\nano_3^{\mu}(\grave{\etab},\grave{a})|
\le C\rhs_{\lip}^{\mu}(\etab,\grave{\etab},a,\grave{a})
\end{equation}
as well.
Combining \eqref{N2 lip}, \eqref{N1 lip}, and \eqref{N3 lip}, we establish \eqref{lip goal}.

\subsubsection{Proof of part \ref{workhorse2 boot}}
Let 
\[
\rhs_{\boot,r}^{\mu}(\etab,a)
:= \mu^{1/2} 
+ \mu^{3/2-r/2}|a| 
+ \norm{\etab}_{r,q_{\star}}
+ \mu^{1-r/2}|a|\norm{\etab}_{r,q_{\star}}
+ \mu^{1-r/2}|a|^2
+ \norm{\etab}_{r,q_{\star}}^2. 
\]
We show that 
\begin{equation}\label{boot goal}
\norm{\nano_1^{\mu}(\etab,a)}_{r+1,q_{\star}}
+ \norm{\nano_2^{\mu}(\etab,a)}_{r+1,q_{\star}}
\le C(r)\rhs_{\boot,r}^{\mu}(\etab,a),
\end{equation}
which will imply \eqref{workhorse2 boot est1}.

First, using the definition of $\nano_2^{\mu}$ in \eqref{N2 N3 eqn defn} and the estimates on $\B_{\mu}$ with $j=1$ from Proposition \ref{A B prop}, we have
\[
\norm{\nano_2^{\mu}(\etab,a)}_{r+1,q_{\star}}
\le \sum_{k=1}^6 \norm{\B_{\mu}\ell_k^{\mu}(\etab,a)}_{r+1,q_{\star}}
\le C(r)\mu^{-1}\sum_{k=1}^6 \norm{\ell_k^{\mu}(\etab,a)}_{r,q_{\star}}.
\]
This, in conjunction with the mapping bounds \eqref{ell map} and \eqref{L map}, gives
\begin{equation}\label{N2 boot}
\norm{\nano_2^{\mu}(\etab,a)}_{r+1,q_{\star}}
\le C(r)\mu^{-1}\Lfrak_{\map,r}^{\mu}(\etab,a)
\le C(r)\rhs_{\boot,r}^{\mu}(\etab,a).
\end{equation}

Next, the definition of $\nano_1^{\mu}$ in \eqref{h5} and the smoothing properties of $\H_c$ imply
\[
\norm{\nano_1^{\mu}(\etab,a)}_{r+1,q_{\star}}
\le \sum_{k=1}^5 \norm{\H_c^{-1}h_k^{\mu}(\etab,a)}_{r+1,q_{\star}}
\le C\sum_{k=1}^4 \norm{h_k^{\mu}(\etab,a)}_{r-1,q_{\star}} + C\norm{\nano_2^{\mu}(\etab,a)}_{r,q_{\star}}.
\]
The mapping bounds \eqref{h map} and \eqref{h map} provide
\[
\norm{h_k^{\mu}(\etab,a)}_{r-1,q_{\star}}
\le C(r)\left(\mu^{1-(r-1)/2}|a| + \mu\norm{\etab}_{r-1,q_{\star}}+\mu^{1-(r-1)/2}|a|\norm{\etab}_{r-1,q_{\star}} + \norm{\etab}_{r-1,q_{\star}}^2\right).
\]
Since $1-(r-1)/2 = 3/2-r/2$, we may further bound
\[
\norm{h_k^{\mu}(\etab,a)}_{r-1,q_{\star}}
\le C(r)\rhs_{\boot,r}^{\mu}(\etab,a),
\]
and thus
\begin{equation}\label{N1 boot}
\norm{\nano_1^{\mu}(\etab,a)}_{r+1,q_{\star}}
\le C(r))\rhs_{\boot,r}^{\mu}(\etab,a).
\end{equation}

By \eqref{N3 map r}, we have
\begin{multline}\label{N3 boot}
|\nano_3^{\mu}(\etab,a)|
\le C(r)\mu^{r/2-2}\Lfrak_{\map,r}^{\mu}(\etab,a) \\
\\
= C(r)\left(\mu^{r/2}
+\mu^3|a| + \mu^{r/2-1}\norm{\etab}_{r,q_{\star}} + \mu^{r/2-1}\norm{\etab}_{r,q_{\star}}^2+|a|\norm{\etab}_{r,q_{\star}} + |a|^2\right). 
\end{multline}
This proves \eqref{workhorse2 boot est2} and completes the proof of Lemma \ref{workhorse lemma for main workhorse lemma}.

\begin{remark}\label{refined leading order limit remark}
Observe from \eqref{L map} that the first term in $\Lfrak_{\map,r}^{\mu}$ is $\mu^2$.
This is an artifact of the refined leading order limit of Theorem \ref{refined leading order limit thm} and appears in the estimate on $\ell_1^{\mu}$ below in Appendix \ref{ell-1-mu ests appendix}.
Had we not made the refined leading order limit, this term $\ell_1^{\mu}$ would be instead $-\G_c(\sigmab_c,\mu)\cdot\e_2$; see the calculations at the start of Section \ref{beale ansatz section} and in particular \eqref{h ell terms}.
As we discussed in Section \ref{refined leading order limit section}, we have only $\norm{\G_c(\sigmab_c,\mu)\cdot\e_2}_{r,q_c} = \O(\mu)$.
But then in taking the estimate \eqref{N2 boot}, multiplying by $\mu^{-1}$ would leave us with an $\O(1)$ term in $\rhs_{\boot,r}^{\mu}(\etab,a)$, and this would ruin the proof of part \ref{main workhorse boot} of Lemma \ref{main workhorse lemma} given in Appendix \ref{proof of main workhorse lemma appendix}.
\end{remark}

\subsection{Estimates on the terms $h_1^{\mu},\ldots,h_5^{\mu}$ and $\ell_1^{\mu},\ldots,\ell_6^{\mu}$}\label{h ell ests appendix}

We recall that the terms $h_1^{\mu},\ldots,h_4^{\mu}$ and $\ell_1^{\mu},\ldots,\ell_3^{\mu}$ were defined in \eqref{h ell terms}, while $\ell_4^{\mu}$ and $\ell_5^{\mu}$ appeared in \eqref{ell-4 ell-5 defns}, $\ell_6^{\mu}$ in \eqref{return to original frequency}, and $h_5^{\mu}$ in \eqref{h5}.

We will rely heavily on the following collection of abstract estimates, which was stated as Proposition G.1 in \cite{faver-hupkes}.
Its proof is a combination of calculus and techniques from \cite{faver-wright}, in particular Lemma A.2 in that paper.

\begin{proposition}\label{abstract sobolev prop}

\begin{enumerate}[label={\bf(\roman*)},ref={(\roman*)}]

\item
If $f \in H_q^r$ and $g \in W^{r,\infty}$, then 
\begin{equation}\label{Hrq-Wrinfty general est}
\norm{fg}_{r,q}
\le \norm{f}_{r,q}\norm{g}_{W^{r,\infty}}.
\end{equation}

\item
If $f$, $g \in H_q^r$, then
\begin{equation}\label{Hrq algebra est}
\norm{fg}_{r,q}
\le \norm{\sech(q\cdot)}_{W^{r,\infty}}\norm{f}_{r,q}\norm{g}_{r,q}.
\end{equation}

\item
If $f \in W^{r,\infty}$ and $\omega \in \R$, then 
\begin{equation}\label{Wrinfty scaling est}
\norm{f(\omega\cdot)}_{W^{r,\infty}}
\le \left(\max_{0 \le k \le r} |\omega|^k\right)\norm{f}_{W^{r,\infty}}.
\end{equation}

\item\label{decay borrowing}
If $f \in H_{q_2}^r$; $g \in W^{r,\infty}$; $\omega$, $\grave{\omega} \in \R$; and $0 < q_1 < q_2$, then
\begin{equation}\label{decay borrowing est}
\norm{f \cdot (g(\omega\cdot)-g(\grave{\omega}\cdot))}_{r,q_1}
\le C(r,q_2-q_1)\left(\max_{0 \le k \le r} \Lip(\partial_x^k[g])\right)\norm{f}_{r,q_2}\norm{g}_{W^{r,\infty}}|\omega-\grave{\omega}|,
\end{equation}
where, for a function $h \colon \R \to \C$, the Lipschitz constant $\Lip(h)$ is
\[
\Lip(h)
:= \sup_{\substack{x,\grave{x} \in \R \\ x \ne \grave{x}}} \left|\frac{h(x)-h(\grave{x})}{x-\grave{x}}\right|.
\]

\item\label{decay borrowing BIG}
If $f \in H_{q_2}^r$; $g$, $\grave{g} \in W^{r,\infty}$; $\omega$, $\grave{\omega} \in \R$; and $0 < q_1 < q_2$, then
\begin{multline}\label{outer inner Lipschitz}
\norm{f \cdot (g(\omega\cdot)-\grave{g}(\grave{\omega}\cdot))}_{r,q_1} \\
\\
\le C(r,q_2-q_1)\norm{f}_{r,q_2}\left[\left(\max_{0 \le k \le r} \Lip(\partial_x^k[g])\right)\norm{g}_{W^{r,\infty}}|\omega-\grave{\omega}|
+ \left(\max_{0 \le k \le r} |\grave{\omega}|^k\right)\norm{g-\grave{g}}_{W^{r,\infty}}\right] .
\end{multline}
\end{enumerate}
\end{proposition}

Using Theorem \ref{periodic solutions thm}, specifically the Lipschitz estimates in \eqref{naked periodic ests}, and the scaling estimate in \eqref{Wrinfty scaling est}, it is straightforward to deduce the following additional estimates on the periodic profiles.

\begin{lemma}\label{periodic cor}

\begin{enumerate}[label={\bf(\roman*)}]

\item
Let $\mu \in (0,\mu_{\per})$.
Then
\begin{equation}\label{actual periodic bound}
\sup_{|a| \le a_{\per}} |\omega_{\mu}^a| + \norm{\phi_{\mu,1}^a}_{W^{r,\infty}} \le C(r)\mu^{1-r/2}.
\end{equation}

\item
Let $\mu \in (0,\mu_{\per})$, $r \ge 0$, and $|a|$, $|\grave{a}| \le a_{\per}$.
Then
\begin{equation}\label{actual periodic lip}
\norm{\phi_{\mu,1}^a-\phi_{\mu,1}^{\grave{a}}}_{W^{r,\infty}} \le C(r)\mu^{1-r/2}|a-\grave{a}|.
\end{equation}

\item
Let $r \ge 0$ and $|a| \le a_{\per}$.
Then
\begin{equation}\label{psi1 a}
\sup_{\mu \in (0,\mu_{\per})} \norm{\psi_{\mu,1}^a}_{W^{r,\infty}} \le C(r)|a|.
\end{equation}

\end{enumerate}
\end{lemma}

Throughout the following estimates on the $h_k^{\mu}$ and $\ell_k^{\mu}$ terms, we recall that $a$, $\grave{a}$, $\mu \in (0,1)$, which will allow us to bound some terms quickly.
Also, since $\qbar_{\star} < q_{\star}$, we have
\[
\norm{\etab}_{r,\qbar_{\star}}
\le \norm{\etab}_{r,q_{\star}}.
\]

\subsubsection{Estimates on $h_1^{\mu}$}
We directly estimate
\[
\norm{h_1^{\mu}(\etab,a)}_{r,q_{\star}}
\le C(r)\mu\norm{\etab}_{r,q_{\star}}
\quadword{and}
\norm{h_1^{\mu}(\etab,a)-h_1^{\mu}(\grave{\etab},\grave{a})}_{3,\qbar_{\star}}
\le C\mu\norm{\etab-\grave{\etab}}_{3,\qbar_{\star}}.
\]

\subsubsection{Estimates on $h_2^{\mu}$}\label{h-2-mu ests appendix}
We use the product estimate \eqref{Hrq-Wrinfty general est}, the scaling estimate \eqref{Wrinfty scaling est}, and the periodic estimate \eqref{actual periodic bound} in Corollary \ref{periodic cor} to bound
\begin{equation}\label{h-2-mu bound}
\norm{h_2^{\mu}(\etab,a)}_{r,q_{\star}}
\le C(r)\mu^{1-r/2}|a|.
\end{equation}

Next, keeping $r$ arbitrary for later convenience, we estimate
\[
\norm{h_2^{\mu}(\etab,a)-h_2^{\mu}(\grave{\etab},\grave{a})}_{r,\qbar_{\star}}
\le \bunderbrace{C|a-\grave{a}|\norm{\varsigma_{c,1}^{\mu}\phi_{\mu,1}^a}_{r,\qbar_{\star}}}{I}
+ \bunderbrace{C|a|\norm{\varsigma_{c,1}^{\mu}\big(\phi_{\mu,1}^a-\phi_{\mu,1}^{\grave{a}}\big)}_{r,\qbar_{\star}}}{II}.
\]
The same estimates that implied \eqref{h-2-mu bound} also give
\[
I
\le C\mu^{1-r/2}|a-\grave{a}|.
\]
For $II$, since $\qbar_{\star} < q_c$, we use the massive decay borrowing estimate \eqref{outer inner Lipschitz} and the various periodic estimates in \eqref{naked periodic ests} and Corollary \ref{periodic cor} to bound
\begin{equation}\label{h-2-mu II lip}
II
\le C\mu^{1-r/2}|\grave{a}||a-\grave{a}|.
\end{equation}

\subsubsection{Estimates on $h_3^{\mu}$}\label{h-3-mu ests appendix}
We use the product estimate \eqref{Hrq-Wrinfty general est} and the periodic estimate \eqref{actual periodic bound} to bound
\begin{equation}\label{h-3-mu bound}
\norm{h_3^{\mu}(\etab,a)}_{r,q_{\star}}
\le C(r)\mu^{1-r/2}|a|\norm{\etab}_{r,q_{\star}}.
\end{equation}

Next, keeping $r$ arbitrary for later convenience, we have
\begin{multline*}
\norm{h_3^{\mu}(\etab,a)-h_3^{\mu}(\grave{\etab},\grave{a})}_{r,\qbar_{\star}}
\le \bunderbrace{C(r,q_{\star})|a-\grave{a}|\norm{\phi_{\mu,1}^a\eta_1}_{r,\qbar_{\star}}}{I}
+ \bunderbrace{|\grave{a}|\norm{\big(\phi_{\mu,1}^a-\phi_{\mu,1}^{\grave{a}}\big)\eta_1}_{r,\qbar_{\star}}}{II} \\
\\
+ \bunderbrace{|\grave{a}|\norm{\phi_{\mu,1}^{\grave{a}}(\eta_1-\grave{\eta}_1)}_{r,\qbar_{\star}}}{III}.
\end{multline*}
As in Appendix \ref{h-2-mu ests appendix}, we first use the decay borrowing estimate \eqref{decay borrowing est}, since $\qbar_{\star} < q_{\star}$, and find
\[
I 
\le C\mu^{1-r/2}|a-\grave{a}|\norm{\etab}_{r,q_{\star}}.
\]
Next, we estimate
\[
III
\le C\mu^{1-r/2}|\grave{a}|\norm{\etab-\grave{\etab}}_{r,\qbar_{\star}}
\]
just as we did with \eqref{h-3-mu bound}.
For $II$, we use decay borrowing and the periodic estimates as we did to establish \eqref{h-2-mu II lip}, and we find
\[
II
\le C\mu^{1-r/2}|\grave{a}|\norm{\etab}_{r,q_{\star}}|a-\grave{a}|.
\]

\subsubsection{Estimates on $h_4^{\mu}$}
Since $h_4^{\mu}$ is quadratic in $\eta_1$, we have
\[
\norm{h_4^{\mu}(\etab,a)}_{r,q_{\star}}
\le C(r)\norm{\etab}_{r,q_{\star}}^2
\]
and
\[
\norm{h_4^{\mu}(\etab,a)-h_4^{\mu}(\grave{\etab},\grave{a})}_{3,\qbar_{\star}}
\le C\norm{\etab+\grave{\etab}}_{r,\qbar_{\star}}\norm{\etab-\grave{\etab}}_{3,\qbar_{\star}}.
\]

\subsubsection{Estimates on $h_5^{\mu}$}
The requisite estimates on $h_5^{\mu}$ are $\O(1)$ multiples of those for the $\ell_k^{\mu}$ terms, so we just proceed to establish those below.

\subsubsection{Estimates on $\ell_1^{\mu}$}\label{ell-1-mu ests appendix}
Theorem \ref{refined leading order limit thm} implies
\[
\norm{\ell_1^{\mu}(\etab,a)}_{r,q_{\star}}
\le C(r,q_{\star})\mu^2,
\]
and since $\ell_1^{\mu}$ is constant in both $\etab$ and $a$, its Lipschitz estimate is trivial.

\subsubsection{Estimates on $\ell_2^{\mu}$}
We directly estimate
\[
\norm{\ell_2^{\mu}(\etab,a)}_{r,q_{\star}}
\le C(r)\mu\left(\norm{\etab}_{r,q_{\star}} + \norm{\etab}_{r,q_{\star}}^2\right)
\]
and
\[
\norm{\ell_2^{\mu}(\etab,a)}_{3,\qbar_{\star}}
\le C\mu\left(\norm{\etab-\grave{\etab}}_{3,\qbar_{\star}} + \norm{\etab+\grave{\etab}}_{3,\qbar_{\star}}\norm{\etab-\grave{\etab}}_{3,\qbar_{\star}}\right).
\]

\subsubsection{Estimates on $\ell_3^{\mu}$}
These estimates are the same as those on $h_3^{\mu}$ from Appendix \ref{h-3-mu ests appendix}, except now we gain an extra factor of $\mu$:
\[
\norm{\ell_3^{\mu}(\etab,a)}_{r,q_{\star}}
\le C(r)\mu^{2-r/2}|a|\norm{\etab}_{r,q_{\star}}
\]
and (taking $r=3$ in the $h_2^{\mu}$ Lipschitz estimates)
\[
\norm{\ell_3^{\mu}(\etab,a)-\ell_3^{\mu}(\grave{\etab},\grave{a})}_{3,\qbar_{\star}}
\le C\mu^{2-r/2}\left(|a-\grave{a}|\norm{\etab}_{3,\qbar_{\star}} + |\grave{a}|\norm{\etab-\grave{\etab}}_{3,\qbar_{\star}} + |\grave{a}||a-\grave{a}|\norm{\etab}_{3,\qbar_{\star}}\right).
\]

\subsubsection{Estimates on $\ell_4^{\mu}$}
We recall from Theorem \ref{periodic solutions thm} that 
\[
\phi_{\mu,1}^a-\phi_{\mu,1}^0
= \mu\upsilon_{\mu}\big(\cos(\omega_{\mu}^a\cdot)-\cos(\omega_{\mu}^0\cdot)\big) + \mu\psi_{\mu,1}^a(\omega_{\mu}^a\cdot),
\]
and so
\[
\ell_4^{\mu}(\etab,a)
= \bunderbrace{2a\mu^2\upsilon_{\mu}\delta\big[\sigma_c\big(\cos(\omega_{\mu}^a\cdot)-\cos(\omega_{\mu}^0\cdot)\big)\big]}{\ell_{41}^{\mu}(\etab,a)}
+ \bunderbrace{2a\mu^2\delta\big[\sigma_c\psi_{\mu,1}^a(\omega_{\mu}^a\cdot)\big]}{\ell_{42}^{\mu}(\etab,a)}.
\]
Since $\upsilon_{\mu} = \O(1)$ and $q_{\star} < q_c$, we may use the decay borrowing estimate \eqref{decay borrowing est} and the periodic Lipschitz estimate \eqref{naked periodic ests} to bound
\[
\norm{\ell_{41}^{\mu}(\etab,a)}_{r,q_{\star}}
\le C(r)\mu^2a^2.
\]
Next, the scaling estimate \eqref{Wrinfty scaling est} and the periodic estimates \eqref{actual periodic bound} and \eqref{psi1 a} imply
\[
\norm{\ell_{42}^{\mu}(\etab,a)}_{r,q_{\star}}
\le C(r,q_{\star})\mu^2|a|\norm{\psi_{\mu,1}^a(\omega_{\mu}^a\cdot)}_{W^{r,\infty}}
\le C(r,q_{\star})\mu^{2-r/2}a^2.
\]
All together, we conclude
\[
\norm{\ell_4^{\mu}(\etab,a)}_{r,q_{\star}}
\le C(r)\mu^{2-r/2}a^2.
\]

For the Lipschitz estimates, we first calculate
\begin{multline*}
\ell_{41}^{\mu}(\etab,a)-\ell_{41}^{\mu}(\grave{\etab},\grave{a})
= \bunderbrace{2\mu^2\upsilon_{\mu}(a-\grave{a})\delta\big[\sigma_c\big(\cos(\omega_{\mu}^a\cdot)-\cos(\omega_{\mu}^0\cdot)\big)\big]}{I} \\
\\
+ \bunderbrace{2\mu^2\upsilon_{\mu}|\grave{a}|\delta\big[\sigma_c\big(\cos(\omega_{\mu}^a\cdot)-\cos(\omega_{\mu}^{\grave{a}}\cdot)\big)\big]}{II}.
\end{multline*}
Then to bound $I$ and $II$, we use the decay borrowing estimate \eqref{decay borrowing est} and the periodic Lipschitz estimate \eqref{naked periodic ests} and find
\[
\norm{I}_{3,q}
\le C\mu^2|a||a-\grave{a}|
\quadword{and}
\norm{II}_{3,q}
\le C\mu^2|\grave{a}||a-\grave{a}|.
\]

\subsubsection{Estimates on $\ell_5^{\mu}$}
These estimates are effectively the same as those on $h_2^{\mu}$ from Appendix \ref{h-2-mu ests appendix}, except now we gain an extra factor of $\mu^2$:
\[
\norm{\ell_5^{\mu}(\etab,a)}_{r,q_{\star}}
\le C(r)\mu^{3-r/2}|a|
\]
and (taking $r=3$ in the $h_2^{\mu}$ Lipschitz estimates)
\[
\norm{\ell_5^{\mu}(\etab,a)-\ell_5^{\mu}(\grave{\etab},\grave{a})}_{3,\qbar_{\star}}
\le C\left(\mu^{3/2}|a-\grave{a}| + \mu^2|\grave{a}|\norm{\etab-\grave{\etab}}_{3,\qbar_{\star}} + \mu^{3/2}|\grave{a}|\norm{\etab}_{r,\qbar_{\star}}|a-\grave{a}|\right).
\]

\subsubsection{Estimates on $\ell_6^{\mu}$}
First, we use the decay borrowing estimate \eqref{decay borrowing est} to find
\begin{multline*}
\norm{\ell_6^{\mu}(\etab,a)}_{r,q_{\star}}
\le C(r,q_{\star})\mu^2|a|\norm{\sigma_c\big(\cos(\omega_{\mu}\cdot)-\cos(\Omega_{\mu}\cdot)\big)}_{r,q_{\star}} \\
\\
\le C(r,q_{\star})\mu^2|a|\norm{\sigma_c}_{r,q_{\star}}|\omega_{\mu}-\Omega_{\mu}|.
\end{multline*}
Next, recall from part \ref{omega vs Omega} of Lemma \ref{omega-mu lemma} that $\omega_{\mu}-\Omega_{\mu} = \O(\mu^{3/2})$, and so we conclude
\[
\norm{\ell_6^{\mu}(\etab,a)}_{r,q_{\star}}
\le C(r,q_{\star})\mu^{7/2}|a|.
\]
The same ideas give the Lipschitz estimate
\[
\norm{\ell_6^{\mu}(\etab,a)-\ell_6^{\mu}(\grave{\etab},\grave{a})}_{3,\qbar_{\star}}
\le C\mu^{7/2}|a-\grave{a}|.
\]

\begingroup
\footnotesize
\setlength{\parskip}{0pt}
\bibliographystyle{siam}
\bibliography{mim-nanopteron-bib}

\begin{thebibliography}{10}

\bibitem{amick-toland}
{\sc C.~J. Amick and J.~F. Toland}, {\em Solitary waves with surface tension.
  {I}. {T}rajectories homoclinic to periodic orbits in four dimensions}, Arch.
  Rational Mech. Anal., 118 (1992), pp.~37--69.

\bibitem{beale1}
{\sc J.~T. Beale}, {\em Water waves generated by a pressure disturbance on a
  steady stream}, Duke Math. J., 47 (1980), pp.~297--323.

\bibitem{beale2}
\leavevmode\vrule height 2pt depth -1.6pt width 23pt, {\em Exact solitary water
  waves with capillary ripples at infinity}, Comm. Pure Appl. Math., 44 (1991),
  pp.~211--257.

\bibitem{bonanomi-et-al}
{\sc L.~Bonanomi, G.~Theocharis, and C.~Daraio}, {\em Wave propagation in
  granular chains with local resonances}, Physical review. E, Statistical,
  nonlinear, and soft matter physics, 91 (2015), p.~033208.

\bibitem{boyd}
{\sc J.~P. Boyd}, {\em Weakly Nonlocal Solitary Waves and Beyond-All-Orders
  Asymptotics}, vol.~442 of Mathematics and Its Applications, Kluwer Academic
  Publishers, Dordrecht, The Netherlands, 1998.

\bibitem{cpkd}
{\sc C.~Chong, M.~A. Porter, P.~G. Kevrekidis, and C.~Daraio}, {\em Nonlinear
  coherent structures in granular crystals}, Journal of Physics: Condensed
  Matter, 29 (2016).

\bibitem{crandall-rabinowitz}
{\sc M.~G. Crandall and P.~H. Rabinowitz}, {\em Bifurcation from simple
  eigenvalues}, J. Functional Analysis, 8 (1971), pp.~321--340.

\bibitem{dauxois}
{\sc T.~Dauxois}, {\em {F}ermi, {P}asta, {U}lam, and a mysterious lady},
  Physics Today, 61 (2008), pp.~55--57.

\bibitem{english-pego}
{\sc J.~M. English and R.~L. Pego}, {\em On the solitary wave pulse in a chain
  of beads}, Proceedings of the American Mathematical Society, 133 (2005),
  pp.~1763--1768.

\bibitem{faver-spring-dimer}
{\sc T.~E. Faver}, {\em Nanopteron-stegoton traveling waves in spring dimer
  {F}ermi-{P}asta-{U}lam-{T}singou lattices}.
\newblock arXiv preprint arXiv:1511.00942. Accepted to Quarterly of Applied
  Mathematics in 2019.

\bibitem{faver-dissertation}
{\sc T.~E. Faver}, {\em Nanopteron-stegoton traveling waves in mass and spring
  dimer {F}ermi-{P}asta-{U}lam-{T}singou lattices}, PhD thesis, Drexel
  University, Philadelphia, PA, May 2018.

\bibitem{faver-goodman-wright}
{\sc T.~E. Faver, R.~H. Goodman, and J.~D. Wright}, {\em Solitary waves in
  mass-in-mass lattices for antiresonance mass ratios.}
\newblock In preparation.

\bibitem{faver-hupkes}
{\sc T.~E. Faver and H.~J. Hupkes}, {\em Micropteron traveling waves in
  diatomic {F}ermi-{P}asta-{U}lam-{T}singou lattices under the equal mass
  limit}.
\newblock Submitted. arXiv preprint arXiv:1909.00770, 2019.

\bibitem{faver-wright}
{\sc T.~E. Faver and J.~D. Wright}, {\em Exact diatomic
  {F}ermi-{P}asta-{U}lam-{T}singou solitary waves with optical band ripples at
  infinity}, SIAM Journal on Mathematical Analysis, 50 (2018), pp.~182--250.

\bibitem{fput-original}
{\sc E.~Fermi, J.~Pasta, and S.~Ulam}, {\em Studies of nonlinear problems},
  Lect. Appl. Math., 12 (1955), pp.~143--56.

\bibitem{friesecke-pego1}
{\sc G.~Friesecke and R.~L. Pego}, {\em Solitary waves on {FPU} lattices. {I}.
  {Q}ualitative properties, renormalization and continuum limit}, Nonlinearity,
  12 (1999), pp.~1601--1627.

\bibitem{friesecke-wattis}
{\sc G.~Friesecke and J.~A.~D. Wattis}, {\em Existence theorem for solitary
  waves on lattices}, Comm. Math. Phys., 161 (1994), pp.~391--418.

\bibitem{kev-vain-et-al}
{\sc P.~G.~Kevrekidis, A.~Vainchtein, M.~Serra-Garcia, and C.~Daraio}, {\em
  Interaction of traveling waves with mass-with-mass defects within a
  {H}ertzian chain}, Physical review. E, Statistical, nonlinear, and soft
  matter physics, 87 (2013), p.~042911.

\bibitem{gmwz}
{\sc J.~Gaison, S.~Moskow, J.~D. Wright, and Q.~Zhang}, {\em Approximation of
  polyatomic {FPU} lattices by {K}d{V} equations}, Multiscale Model. Simul., 12
  (2014), pp.~953--995.

\bibitem{gantzounis}
{\sc G.~Gantzounis, M.~Serra-Garcia, K.~Homma, J.~M. Mendoza, and C.~Daraio},
  {\em Granular metamaterials for vibration mitigation}, Journal of Applied
  Physics, 114 (2013), p.~093514.

\bibitem{co-ops}
{\sc N.~Giardetti, A.~Shapiro, S.~Windle, and J.~D. Wright}, {\em Metastability
  of solitary waves in diatomic {FPUT} lattices}, Mathematics in Engineering, 1
  (2019), pp.~419--433.

\bibitem{herrmann-matthies-asymptotic}
{\sc M.~Herrmann and K.~Matthies}, {\em Asymptotic formulas for solitary waves
  in the high-energy limit of {F}{P}{U}-type chains}, Nonlinearity, 28 (2015),
  pp.~2767--2789.

\bibitem{herrmann-matthies-uniqueness}
\leavevmode\vrule height 2pt depth -1.6pt width 23pt, {\em Uniqueness of
  solitary waves in the high-energy limit of {FPU}-type chains}, in Patterns of
  dynamics, G.~P., H.~J., S.~B., and S.~A., eds., vol.~205 of Springer
  Proceedings in Mathematics \& Statistics, Springer, 2017, pp.~3--15.

\bibitem{herrmann-matthies-stability}
\leavevmode\vrule height 2pt depth -1.6pt width 23pt, {\em Stability of
  high-energy solitary waves in {F}ermi-{P}asta-{U}lam-{T}singou chains},
  Trans. Amer. Math. Soc., 372 (2019), pp.~3425--3486.

\bibitem{hoffman-wright}
{\sc A.~Hoffman and J.~D. Wright}, {\em Nanopteron solutions of diatomic
  {F}ermi-{P}asta-{U}lam-{T}singou lattices with small mass-ratio}, Physica D:
  Nonlinear Phenomena, 358 (2017), pp.~33--59.

\bibitem{iooss}
{\sc G.~Iooss}, {\em Travelling waves in the {F}ermi-{P}asta-{U}lam lattice},
  Nonlinearity, 13 (2000), pp.~849--866.

\bibitem{IK}
{\sc G.~Iooss and K.~Kirchg\"{a}ssner}, {\em Travelling waves in a chain of
  coupled nonlinear oscillators}, Commun. Math, Phys, 211 (2000), pp.~439--464.

\bibitem{ji-hong}
{\sc J.-Y. Ji and J.~Hong}, {\em Existence criterion of solitary waves in a
  chain of grains}, Physics Letters A, 260 (1999), pp.~60--61.

\bibitem{johnson-wright}
{\sc M.~A. Johnson and J.~D. Wright}, {\em Generalized solitary waves in the
  gravity-capillary {W}hitham equation}.
\newblock arXiv preprint arXiv:1807.1149v1.

\bibitem{kevrekidis}
{\sc P.~G. Kevrekidis}, {\em Non-linear waves in lattices: past, present,
  future}, IMA J. Appl. Math., 76 (2011), pp.~389--423.

\bibitem{ksx}
{\sc P.~G. Kevrekidis, A.~G. Stefanov, and H.~Xu}, {\em Traveling waves for the
  mass in mass model of granular chains}, Letters in Mathematical Physics, 106
  (2016), pp.~1067--1088.

\bibitem{woodpile}
{\sc E.~Kim, F.~Li, C.~Chong, G.~Theocharis, J.~Yang, and P.~G. Kevrekidis},
  {\em Highly nonlinear wave propagation in elastic woodpile periodic
  structures}, Phys. Rev. Lett., 114 (2015), p.~118002.

\bibitem{Li}
{\sc F.~Li, P.~Anzel, J.~Yang, P.~G.~Kevrekidis, and C.~Daraio}, {\em Granular
  acoustic switches and logic elements}, Nature communications, 5 (2014),
  p.~5311.

\bibitem{lombardi}
{\sc E.~Lombardi}, {\em Oscillatory Integrals and Phenomena Beyond all
  Algebraic Orders with Applications to Homoclinic Orbits in Reversible
  Systems}, vol.~1741 of Lecture Notes in Mathematics, Springer-Verlag Berlin
  Heidelberg, 2000.

\bibitem{lustri-porter}
{\sc C.~Lustri and M.~A. Porter}, {\em Nanoptera in a period-2 {T}oda chain},
  SIAM J. Applied Dynamical Systems, 17 (2018), pp.~1182--1212.

\bibitem{lustri}
{\sc C.~J. Lustri}, {\em Nanoptera and {S}tokes curves in the 2-periodic
  {F}ermi-{P}asta-{U}lam-{T}singou equation}.
\newblock arXiv preprint arXiv:1905.07092.

\bibitem{mackay}
{\sc R.~S. MacKay}, {\em Solitary waves in a chain of beads under {H}ertz
  contact}, Physics Letters A, 251 (1999), pp.~191--192.

\bibitem{nesterenko-propagation}
{\sc V.~F. Nesterenko}, {\em Propagation of nonlinear compression pulses in
  granular media}, Journal of Applied Mechanics and Technical Physics, 24
  (1983), pp.~733--743.

\bibitem{nesterenko}
\leavevmode\vrule height 2pt depth -1.6pt width 23pt, {\em Dynamics of
  Heterogeneous Materials}, Shock Wave and High Pressure Phenomena, Springer,
  2001.

\bibitem{cement}
{\sc X.~Ni, P.~Rizzo, J.~Yang, D.~Katri, and C.~Daraio}, {\em Monitoring the
  hydration of cement using highly nonlinear solitary waves}, NDT and E
  International, 52 (2012).

\bibitem{pankov}
{\sc A.~Pankov}, {\em Travelling Waves and Periodic Oscillations in
  Fermi-Pasta-Ulam Lattices}, Imperial College Press, Singapore, 2005.

\bibitem{Spadoni}
{\sc A.~Spadoni and C.~Daraio}, {\em Generation and control of sound bullets
  with a nonlinear acoustic lens}, Proceedings of the National Academy of
  Sciences, 107 (2010), pp.~7230--7234.

\bibitem{sv}
{\sc Y.~Starosvetsky and A.~Vainchtein}, {\em Solitary waves in {F}{P}{U}
  lattices with alternating bond potentials}, Mechanics Research
  Communications,  (2017).

\bibitem{sk12}
{\sc A.~Stefanov and P.~Kevrekidis}, {\em On the existence of solitary
  traveling waves for generalized {H}ertzian chains}, J. Nonlinear Sci., 22
  (2012), pp.~327--349.

\bibitem{sk13}
\leavevmode\vrule height 2pt depth -1.6pt width 23pt, {\em Traveling waves for
  monomer chains with precompression}, Nonlinearity, 26 (2013), pp.~539--564.

\bibitem{sun-shen}
{\sc S.~M. Sun and M.~C. Shen}, {\em Exponentially small estimate for the
  amplitude of capillary ripples of a generalized solitary wave}, J. Math.
  Anal. Appl., 172 (1993), pp.~533--566.

\bibitem{toda}
{\sc M.~Toda}, {\em Theory of nonlinear lattices}, vol.~20 of Springer Series
  in Solid State Sciences, Springer, Berlin Heidelberg, 1978.

\bibitem{vswp}
{\sc A.~Vainchtein, Y.~Starosvetsky, J.~D. Wright, and R.~Perline}, {\em
  Solitary waves in diatomic chains}, Phys. Rev. E, 93 (2016).

\bibitem{venney-zimmer}
{\sc C.~R. Venney and J.~Zimmer}, {\em Travelling lattice waves in a toy model
  of {L}ennard-{J}ones interaction}, Quarterly of Applied Mathematics, 72
  (2014), pp.~65--84.

\bibitem{vorotnikov-et-al}
{\sc K.~Vorotnikov, Y.~Starosvetsky, G.~Theocharis, and P.~G. Kevrekidis}, {\em
  Wave propagation in a strongly nonlinear locally resonant granular crystal},
  Physica D: Nonlinear Phenomena, 365 (2018), pp.~27--41.

\bibitem{xu}
{\sc H.~Xu, P.~G.~Kevrekidis, and A.~Stefanov}, {\em Traveling waves and their
  tails in locally resonant granular systems}, Journal of Physics A:
  Mathematical and Theoretical, 48 (2015).

\bibitem{bones}
{\sc J.~Yang, S.~N. Sangiorgio, S.~L. Borkowski, C.~Silvestro, L.~De~Nardo,
  C.~Daraio, and E.~Ebramzadeh}, {\em Site-specific quantification of bone
  quality using highly nonlinear solitary waves}, J Biomech Eng., 134 (2012),
  pp.~101001--101001--8.

\bibitem{zeidler}
{\sc E.~Zeidler}, {\em Applied functional analysis}, vol.~109 of Applied
  Mathematical Sciences, Springer-Verlag, New York, 1995.
\newblock Main principles and their applications.

\end{thebibliography}
\endgroup
\end{document}